\documentclass[10pt, leqno]{article}
\usepackage{amsmath,amsfonts,amsthm,amscd,amssymb}
\numberwithin{equation}{section}
\theoremstyle{plain}
\newtheorem{theo}{Theorem}[section]
\newtheorem{prop}[theo]{Proposition}
\newtheorem{coro}[theo]{Corollary} 
\newtheorem{lemm}[theo]{Lemma}
\theoremstyle{definition}
\newtheorem{defi}[theo]{Definition}
\newtheorem{rema}[theo]{Remark}
\newtheorem{theo-defi}[theo]{Theorem-Definition}
\newtheorem{prop-defi}[theo]{Proposition-Definition}
\newtheorem{rema-defi}[theo]{Remark-Definition}
\newtheorem{exem-defi} [theo]{Example-Definiton}

\def \al{\alpha}

\def \bul{\bullet}
\def \col{\colon}
\def \Del{\Delta}

\def \eps{\epsilon}
\def \Gam{\Gamma}

\def \inf{\infty}

\def \Lam{\Lambda}
\def \lam{\lambda}

\def \lo{\longrightarrow}
\def \lom{\longmapsto}
\def \mab{\mathbb}

\def \Om{\Omega}
\def \om{\omega}
\def \ol{\overline}
\def \os{\overset}
\def \parno{\par\noindent}

\def \sus{\subset}

\def \ul{\underline}

\def \wt{\widetilde}
\bigskip

\newcommand{\getsfrom}{\ensuremath{
\longleftarrow\kern-.50em\lower.0ex\hbox%
{$\shortmid\,$}}}

\begin{document}
\title{An ideal proof for Fujisawa's result and its generalization}
\author{Yukiyoshi Nakkajima
\date{}\thanks{\endgraf}}
\maketitle

\bigskip
\parno
{\bf Abstract.---}
We give a generalization of Fujisawa's theorem in \cite{fut} in a stronger form. 
Our proof of the generalized theorem is purely algebraic 
and it is simpler than his proof. 

\section{Introduction}
This article is a sequel of my article \cite{nhi} in a different subject. 
The principal aim in this article is 
to give an additional result to Fujisawa's article \cite{fut} by giving an ideal proof for it in a stronger form. 
The proof for this result is useful 
as an ideal proof for a part of his main theorem in [loc.~cit.]. 
His proof is not algebraic; 
our proof is purely algebraic and much simpler than his proof in [loc.~cit.]. 
Because of the algebraicity of our proof, we can give a generalization of 
his result without difficulty.  
First let us recall his theorem. 
\par 
Let $r$ be a fixed positive integer. 
Let $s_{\mab C}$ be a log analytic space whose underlying analytic space is 
${\rm Spec}({\mab C})^{\rm an}$ and 
whose log structure is associated to a morphism 
${\mab N}^r\owns e_i\lom 0\in {\mab C}$ of monoids,  
where $e_i$ $(1\leq i\leq r)$ is a canonical basis of  ${\mab N}^r$. 
Let $M_{s_{\mab C}}$ be the log structure of $s_{\mab C}$. 
Set $\ol{M}_{s_{\mab C}}:=M_{s_{\mab C}}/{\mab C}^*$.  
Let $t_i$ be the section of $\ol{M}_{s_{\mab C}}$ corresponding to $e_i$. 
Let $\wt{t}_i$ be a section of $M_{s_{\mab C}}$ which is a lift of $t_i$. 
Consider a section $d\log t_i:=d\log \wt{t}_i\in \Om^1_{s_{\mab C}/{\mab C}}$.  
(This indeed depends only on $t_i$.) 
Let $f\col X\lo s_{\mab C}$ be a (not necessarily proper) 
log smooth morphism of reduced log analytic spaces. 
Let $\os{\circ}{X}$ be the underlying analytic space of $X$. 
Assume that there exist local generators of the ideal sheaves of  
the irreducible components of $\os{\circ}{X}$ 
which are the images of local sections 
of the log structure $M_X$ of $X$ by the structural morphism 
$\al_X\col M_X\lo {\cal O}_X$. 
Let $\os{\circ}{\Del}:=\{x\in {\mab C}~\vert~\vert x \vert <1\}$ be the unit disk. 
Assume that any irreducible component
$\os{\circ}{X}_{\lam}$ $(\lam \in \Lam)$ of $\os{\circ}{X}$ 
is smooth over ${\mab C}$ and that $\os{\circ}{X}$ fits locally into 
the following cartesian diagram 
\begin{equation*} 
\begin{CD}
\os{\circ}{X}@>{\subset}>> \os{\circ}{\Del}{}^n\\
@VVV @VVV \\
\os{\circ}{s_{\mab C}} @>{\subset}>>  \os{\circ}{\Del}{}^r
\end{CD}
\tag{1.0.1}\label{cd:kdfsx}
\end{equation*} 
such that any $\os{\circ}{X}_{\lam}$ in this diagram 
is defined by equations $x_{i_1}=\cdots =x_{i_k}=0$ 
for some $1\leq i_1<\cdots <i_k\leq n$, where $x_1,\ldots, x_n$ are 
the standard coordinates of $\os{\circ}{\Del}{}^n$. 
If $X$ is locally a finitely many product of 
SNCL(=simple normal crossing log) analytic spaces 
over the log point of ${\rm Spec}({\mab C})^{\rm an}$, 
then these assumptions are satisfied. 
\par 
Let $\Om^{\bul}_{X/s_{\mab C}}$ be the log de Rham complex of $X/s_{\mab C}$. 
Let $F$ be the log Hodge filtration on $\Om^{\bul}_{X/s_{\mab C}}$: 
$F^i\Om^{\bul}_{X/s_{\mab C}}:=\Om^{\bul \geq i}_{X/s_{\mab C}}$. 
Let $\os{\circ}{X}{}^{(m)}$ $(m\in {\mab N})$ be the disjoint union of $(m+1)$-fold intersections of the different irreducible components of $\os{\circ}{X}$. 
Endow $\os{\circ}{X}{}^{(m)}$ with the inverse image of the log structure 
of $X$ by the natural morphism $a^{(m)}\col \os{\circ}{X}{}^{(m)}\lo \os{\circ}{X}$ 
and let $X^{(m)}$ be the resulting log analytic space. 
This space is a log analytic space over $s_{\mab C}$ and $\os{\circ}{s}_{\mab C}$. 
Let $f^{(m)}\col X^{(m)}\lo s_{\mab C}$ be the structural morphism. 
Let $\Om^{\bul}_{X^{(m)}/{\mab C}}$ be the log de Rham complex of 
$X^{(m)}/\os{\circ}{s}_{\mab C}$. 
By abuse of notation, denote by $d\log t_i\in 
\Om^1_{X^{(m)}/{\mab C}}$ the image of 
$d\log t_i$ in $\Om^1_{X^{(m)}/{\mab C}}$ by the morphism 
$f^{(m)*}\col f^{(m)*}(\Om^1_{s_{\mab C}/{\mab C}})\lo \Om^1_{X^{(m)}/{\mab C}}$.  
Let $\Om^{\bul}_{X^{(m)}/{\mab C}}[u_1,\ldots,u_r]$ 
be the Hirsch extension of $\Om^{\bul}_{X^{(m)}/{\mab C}}$ by a
morphism ${\mab C}^r\owns u_i\lom d\log t_i\in \Om^1_{X^{(m)}/{\mab C}}$, 
where $u_i:=(0,\ldots,0,\os{i}{1},0,\ldots,0)$. 
(The notion of the Hirsch extension of a dga by a vector space appears in the definition of 
a minimal model of a dga in Sullivan's theory (\cite{su}).) 
That is, $\Om^{\bul}_{X^{(m)}/{\mab C}}[u_1,\ldots,u_r]$ 
is a complex of ${\mab C}$-vector spaces on $\os{\circ}{X}{}^{(m)}$ defined by 
$\Om^i_{X^{(m)}/{\mab C}}[u_1,\ldots,u_r]:=
{\mab C}[u_1,\ldots,u_r]\otimes_{\mab C}\Om^i_{X^{(m)}/{\mab C}}$ 
with differential morphism induced by the derivative 
$d\col \Om^i_{X^{(m)}/{\mab C}}\lo \Om^{i+1}_{X^{(m)}/{\mab C}}$ 
and $u_i^{e}\otimes 1\lom eu_i^{e-1}\otimes d\log t_i$ $(e\in {\mab N})$. 
By abuse of notation, denote also by $d$ the differential morphism 
$\Om^i_{X^{(m)}/{\mab C}}[u_1,\ldots,u_r]\lo 
\Om^{i+1}_{X^{(m)}/{\mab C}}[u_1,\ldots,u_r]$. 
(In our mind, we consider $u_i$ as ``$\log \wt{t}_i$'' only in this introduction 
which can be considered as a function of the universal cover of $\Del \setminus \{O\}$.) 
Fix a total order on $\Lam$. 
Let $H(X/s_{\mab C})$ be the single complex of 
the double complex 
$\{a^{(m)}_*(\Om^i_{X^{(m)}/{\mab C}}[u_1,\ldots,u_r])\}$ $(m,i\in {\mab N})$  
with the following horizontal and vertical differential morphisms as in \cite{nh3}: 
\begin{equation*} 
\begin{CD}
@. \cdots @>>> \cdots\\
@. @AAA @AAA\\
\cdots @>{\rho^{(m-1)*}}>>a^{(m)}_*(\Om^{i+1}_{X^{(m)}/{\mab C}}[u_1,\ldots,u_r])@>{\rho^{(m)*}}>>
a^{(m+1)}_*(\Om^{i+1}_{X^{(m+1)}/{\mab C}}[u_1,\ldots,u_r]) @>{\rho^{(m+1)*}}>>\cdots\\
@. @A{(-1)^md}AA @AA{(-1)^{m+1}d}A \\
\cdots @>{\rho^{(m-1)*}}>>
a^{(m)}_*(\Om^i_{X^{(m)}/{\mab C}}[u_1,\ldots,u_r])@>{\rho^{(m)*}}>>
a^{(m+1)}_*(\Om^{i}_{X^{(m+1)}/{\mab C}}[u_1,\ldots,u_r]) @>{\rho^{(m+1)*}}>>\cdots\\
@. @AAA @AAA\\
@. \cdots @>>> \cdots,\\
\end{CD}
\end{equation*} 
where $\rho^{(m)*}$ is the standard \v{C}ech morphism which will be recalled in the text.  
Set $\ul{u}^{\ul{e}}:=u_1^{e_1}\cdots u_r^{e_r}$ and 
$\vert \ul{e}\vert=\sum_{i=1}^re_i$, 
where $\ul{e}:=(e_1,\ldots,e_r)$. 
(The complex $H(X/s_{\mab C})$ has been denoted by $K_{\mab C}$ 
in \cite{fut}.) 
Set 
\begin{align*} 
F^i(\Om^{\bul}_{X^{(m)}/{\mab C}}[u_1,\ldots,u_r])
:=\bigoplus_{\ul{e}\in {\mab N}^r}{\mab C}\ul{u}^{\ul{e}}\otimes_{\mab C} 
\Om^{\bul \geq i-\vert \ul{e}\vert}_{X^{(m)}/{\mab C}}.
\end{align*}
This filtration indeed induces a filtration $F$ on $H(X/s_{\mab C})$. 
Let $\Om^{\bul}_{X^{(m)}/s_{\mab C}}$ be the log de Rham complex of 
$X^{(m)}/s_{\mab C}$ and let $s(a^{(\star)}_*(\Om^{\bul}_{X^{(\star)}/s_{\mab C}}))$ 
be the single complex of the double complex 
$\{a^{(m)}_*(\Om^i_{X^{(m)}/s_{\mab C}})\}_{m,i\in {\mab N}}$. 
The complex $s(\Om^{\bul}_{X^{(\star)}/s_{\mab C}})$ 
has a natural log Hodge filtration $F$: $F^is(\Om^{\bul}_{X^{(\star)}/s_{\mab C}})
:=s(\Om^{\bul \geq i}_{X^{(\star)}/s_{\mab C}})$.  
Let 
$(H(X/s_{\mab C}),F)\lo (s(a^{(\star)}_*(\Om^{\bul}_{X^{(\star)}/s_{\mab C}})),F)$ 
be the natural morphism defined by the morphism $u_i\lom 0$ and 
the projection 
$\Om^{\bul}_{X^{(\star)}/{\mab C}}\lo 
\Om^{\bul}_{X^{(\star)}/s_{\mab C}}$. 
In his article \cite{fut} Fujisawa has proved that 
the natural morphism 
\begin{align*} 
(\Om^{\bul}_{X/s_{\mab C}},F)\lo (s(a^{(\star)}_*(\Om^{\bul}_{X^{(\star)}/s_{\mab C}})),F)
\tag{1.0.2}\label{ali:korsx}
\end{align*} 
is a filtered quasi-isomorphism. 
Consequently he has proved that there exists a filtered morphism 
\begin{align*} 
(H(X/s_{\mab C}),F)\lo (\Om^{\bul}_{X/s_{\mab C}},F)
\tag{1.0.3}\label{ali:kcfsx}
\end{align*} 
in the derived category ${\rm D}^+{\rm F}({\mab C}_X)$ of
bounded below filtered complexes of sheaves of ${\mab C}$-vector spaces on 
$\os{\circ}{X}$. (Strictly speaking, this has not been stated in [loc.~cit.]
because the category  ${\rm D}^+{\rm F}({\mab C}_X)$ has not appeared in [loc.~cit.]
(cf.~\cite[Corollary 5.7]{fut}).) 
Here note that the filtration $F$ on $H(X/s_{\mab C})$ is not finite;  
we need a formalism of derived categories of filtered complexes in e.~g., \cite{nh2}, 
which is a special case of the formalism of derived categories of 
filtered complexes in \cite{blec}. 
Furthermore he has proved that the morphism (\ref{ali:kcfsx}) induces an isomorphism 
\begin{align*} 
{\rm gr}^i_FH(X/s_{\mab C})\os{\sim}{\lo} {\rm gr}^i_F\Om^{\bul}_{X/s_{\mab C}}. 
\tag{1.0.4}\label{ali:kcgrsx}
\end{align*}  
By this fact, we see that there exists a large number $N$ such that 
there exists an isomorphism 
\begin{align*} 
H(X/s_{\mab C})/F^NH(X/s_{\mab C}) \os{\sim}{\lo} \Om^{\bul}_{X/s_{\mab C}}
\tag{1.0.5}\label{ali:kckffsx}
\end{align*} 
if $\os{\circ}{X}$ is {\it proper} over $\os{\circ}{s}_{\mab C}$. 
We think that the subsheaf $F^NH(X/s_{\mab C})$ should be 
removed if possible. 
By using the isomorphism (\ref{ali:kckffsx}) 
and by constructing a ${\mab Q}$-structure $H_{\mab Q}$ 
of $H(X/s_{\mab C})$ with a weight filtration $W$ and by using 
his generalization on mixed Hodge complexes and 
a fundamental theorem on mixed Hodge structures, 
he has proved the following: 

\begin{theo}[{\rm {\bf \cite[p.~94]{fut}}}]\label{thoe:fj}
Assume that the irreducible components of $\os{\circ}{X}$ are proper and K\"{a}hler. 
Then the morphism {\rm (\ref{ali:kckffsx})}  
induces the following isomorphism of cohomologies with Hodge filtrations: 
\begin{align*} 
(H^q(X,H(X/s_{\mab C})),F)
\os{\sim}{\lo} (H^q(X,\Om^{\bul}_{X/s_{\mab C}}),F).
\tag{1.1.1}\label{ali:kcfcsx}
\end{align*}   
\end{theo}
\par 
In this article we prove that the filtered morphism 
\begin{align*} 
(H(X/s_{\mab C}),F)\lo (\Om^{\bul}_{X/s_{\mab C}},F)
\tag{1.1.2}\label{ali:kcfosx}
\end{align*} 
is indeed a filtered isomorphism in the filtered 
derived category ${\rm D}^+{\rm F}({\mab C}_X)$ 
without assuming that the irreducible components of 
$\os{\circ}{X}$ are proper nor K\"{a}hler 
(we assume only that the irreducible components of $\os{\circ}{X}$ are smooth) 
and without using any result on generalized mixed Hodge complex in \cite{fut}.
Consequently we give a much simpler and algebraic proof of the isomorphism 
(\ref{ali:kcfcsx}) without the assumption in (\ref{thoe:fj}). 
In fact, we give the isomorphism at the level of the filtered complexes and  
we give a generalization of this result for certain coefficients as follows:

\begin{theo}\label{theo:lsl}
Let $X/s_{\mab C}$ be as before {\rm (\ref{thoe:fj})}. 
Let ${\cal E}$ be a locally free coherent ${\cal O}_X$-module and 
let $\nabla  \col {\cal E}\lo {\cal E}\otimes_{{\cal O}_X}\Om^1_{X/{\mab C}}$ be 
a locally nilpotent integrable connection on $X$ with respect to $s_{\mab C}$ 
$($see the text for the definition of this notion$)$. 
Let $\nabla \col {\cal E}\lo {\cal E}\otimes_{{\cal O}_X}\Om^1_{X/s_{\mab C}}$ 
be the induced connection. 
Endow ${\cal E}\otimes_{{\cal O}_X}\Om^{\bul}_{X/s_{\mab C}}$ 
with a filtration $F$ on ${\cal E}\otimes_{{\cal O}_X}\Om^{\bul}_{X/s_{\mab C}}$ 
as follows$:$
$F^i({\cal E}\otimes_{{\cal O}_X}\Om^{\bul}_{X/s_{\mab C}}):=
{\cal E}\otimes_{{\cal O}_X}\Om^{\bul \geq i}_{X/s_{\mab C}}$.  
Set 
$$H(X/s_{\mab C},{\cal E}):=
s({\cal E}\otimes_{{\cal O}_X}a^{(\star)}_*
(\Om^{\bul}_{X^{(\star)}/{\mab C}}[u_1,\ldots,u_r]))$$  
with already essentially defined differential morphism 
and 
\begin{align*} 
F^i({\cal E}\otimes_{{\cal O}_X}
\Om^{\bul}_{X^{(\star)}/{\mab C}}[u_1,\ldots,u_r])
:=\bigoplus_{\ul{e}\in {\mab N}^r}{\mab C}\ul{u}^{\ul{e}}\otimes_{\mab C} 
{\cal E}\otimes_{{\cal O}_X}\Om^{\bul \geq i-\vert \ul{e}\vert}_{X^{(\star)}/{\mab C}}.
\end{align*}
Then there exists the following filtered isomorphism 
\begin{align*} 
(H(X/s_{\mab C},{\cal E}),F)
\os{\sim}{\lo}  ({\cal E}\otimes_{{\cal O}_X}\Om^{\bul}_{X/s_{\mab C}},F). 
\tag{1.2.1}\label{ali:yscc}
\end{align*} 
Moreover this isomorphism is contravariantly functorial 
with respect to certain morphisms of $X/s_{\mab C}$'s and $({\cal E},\nabla)$'s. 
\end{theo}

\parno
This theorem is a special case of the main result (\ref{theo:mra}) below.
As a corollary we obtain the following without assuming 
that $\os{\circ}{X}$ is proper 
over ${\mab C}$ nor without assuming that the irreducible components of 
$\os{\circ}{X}$ are K\"{a}hler: 

\begin{coro}\label{coro:nk}
There exists the following contravariantly functorial isomorphism  
of filtered complexes with respect respect to certain morphisms of 
$X/s_{\mab C}$'s and $({\cal E},\nabla)$'s$:$
\begin{align*} 
R\Gam(X,(H(X/s_{\mab C},{\cal E}),F))
\os{\sim}{\lo} 
R\Gam(X,({\cal E}\otimes_{{\cal O}_X}\Om^{\bul}_{X/s_{\mab C}}),F)).  
\tag{1.3.1}\label{ali:ysoc}
\end{align*} 
Consequently there exists the following isomorphism 
\begin{align*} 
(H^q(X,H(X/s_{\mab C},{\cal E})),F)
\os{\sim}{\lo} 
(H^q(X,{\cal E}\otimes_{{\cal O}_X}\Om^{\bul}_{X/s_{\mab C}}),F)\quad (q\in {\mab Z}).  
\tag{1.3.2}\label{ali:ychsoc}
\end{align*} 
\end{coro}
(\ref{ali:ychsoc}) is a generalization of (\ref{ali:kcfcsx}).
The integrable connection $({\cal E},\nabla)$ is obtained by a local system on $\os{\circ}{X}$. 
More generally ${\cal E}$ is obtained by a locally unipotent local system on $X^{\log}$ 
by the log Riemann-Hilbert correspondence established in \cite{kn}. 
\par
To prove the theorem (\ref{theo:lsl}), we use theory of (PD-)Hirsch extensions   
in \cite{nhi} and a key theorem obtained in [loc.~cit.]. 
\par 
See \cite{nhi} for the log crystalline analogue of (\ref{theo:lsl}) for the case $r=1$. 
In [loc.~cit] we have had to prove the log crystalline analogue of the isomorphism 
\begin{align*} 
H(X/s_{\mab C},{\cal E})
\os{\sim}{\lo}  {\cal E}\otimes_{{\cal O}_X}\Om^{\bul}_{X/s_{\mab C}} 
\tag{1.3.4}\label{ali:yscacc}
\end{align*} 
for the case $r=1$ in order to prove the compatibility of the weight filtration on 
the log crystalline cohomology sheaf of a proper SNCL scheme 
in characteristic $p>0$ with the cup product of the log crystalline cohomology sheaf 
when $({\cal E},\nabla)$ has no log poles. 
Because there is no good analogue of the Hodge filtration on this cohomology sheaf, 
we have had to give a purely algebraic proof for 
the log crystalline analogue of the isomorphism (\ref{ali:yscacc}). 
\par 
The contents of this article are as follows.
\par 
In \S\ref{sec:hi}
we recall results on (PD-)Hirsch extensions in \cite{nhi}. 
\par 
In \S\ref{sec:da} we recall a result in \cite{nhi} and 
we give a key theorem in this article. 
It is an analytic version of a generalization of 
a result in \cite{nhi}, which plays a starting theorem in this article. 
\par 
In \S\ref{sec:mr} we give the main result in this article, 
which is a generalization of (\ref{theo:lsl}).

\section{Recall on the Hirsch extension}\label{sec:hi}
In this section we recall results on the (PD-)Hirsch extension in \cite{nhi}. 
\par 
Let $n$ be a positive integer.  
Let $A$ be a commutative ring with unit element. 
Let $M$ be a projective $A$-module.  
Let $\Gam_A(M)$ be the PD-algebra generated by $M$. 
Let $B_i$ $(1\leq i\leq n)$ be an $A$-algebra.  
Let $\Om^{\bul}_i$  be a dga over $A$ 
such that each $\Om^q_i$ $(q\in {\mab Z})$ is a $B_i$-module.  
Let $E_i$ be a $B_i$-module. 
Let $\varphi_i \col M\lo {\rm Ker}(d\col \Om^1_i\lo \Om^2_i)$ 
be a morphism of $A$-modules. 
Let $(\bigoplus_{i=1}^n(E_i\otimes_{B_i}\Om^{\bul}_i),d)$ 
be a complex of $A$-modules. 
Here $d$ is not necessarily the direct sum of certain boundary morphisms of 
$E_i\otimes_{B_i}\Om^{\bul}_i$.
Then the following $A$-linear morphism 
$$d_H:=d_{H,(\varphi_1,\ldots,\varphi_r)}\col \Gam_A(M)\otimes_A
(\bigoplus_{i=1}^n(E_i\otimes_{B_i}\Om^q_i))
\lo \Gam_A(M)\otimes_A(\bigoplus_{i=1}^n(E_i\otimes_{B_i}\Om^{q+1}_i))$$ 
defined by the following formula 
\begin{align*} 
&d_H(\sum_{i=1}^n\sum_{j(i)_1,\ldots, j(i)_r\in {\mab N}}m^{[j(i)_1]}_1
\cdots m^{[j(i)_r]}_r\otimes e_i\otimes \om_i)\\
&:=
\sum_{i=1}^n\sum_{j(i)_1,\ldots, j(i)_r\in {\mab N}}\sum_{k=1}^r
m^{[j(i)_1]}_1\cdots m^{[j(i)_k-1]}_k \cdots m^{[j(i)_r]}_r\otimes e_i\otimes 
\varphi_M(m_k)\wedge \om_i\\
&+\sum_{i=1}^n\sum_{j(i)_1,\ldots, j(i)_r\in {\mab N}}
m^{[j(i)_1]}_1\cdots m^{[j(i)_r]}_r\otimes d(e_i\otimes \om_i) \\
&\quad (m_1,\ldots m_r\in M,e_i\in E_i,\om_i \in \Om^q_i)
\end{align*} 
makes $\Gam_A(M)\otimes_A(\bigoplus_{i=1}^n(E_i\otimes_{B_i}\Om^{\bul}_i))$ 
a complex of $A$-modules. 
We call the following natural injective morphism 
\begin{align*} 
\bigoplus_{i=1}^nE_i\otimes_{B_i}\Om^{\bul}_i
\os{\sus}{\lo} \Gam_A(M)\otimes_A(\bigoplus_{i=1}^n(E_i\otimes_{B_i}\Om^{\bul}_i))
\end{align*} 
of $A$-modules the {\it PD-Hirsch extension} of 
$(\bigoplus_{i=1}^n(E_i\otimes_{B_i}\Om^{\bul}_i),d)$ by $(M,(\varphi_1,\ldots,\varphi_r))$. 
We denote it by $(\bigoplus_{i=1}^n(E_i\otimes_{B_i}\Om^{\bul}_i))\langle M\rangle
=(\bigoplus_{i=1}^n(E_i\otimes_{B_i}\Om^{\bul}_i))\langle M
\rangle_{(\varphi_1,\ldots,\varphi_r)}$. 
\par 
Let $A\lo A'$ and $B\lo B'$ be morphisms of commutative rings with unit elements. 
Assume that $B'$ is an $A'$-algebra and that the following diagram 
\begin{equation*} 
\begin{CD} 
A@>>> A'\\
@VVV @VVV\\
B@>>>B'
\end{CD}
\end{equation*} 
is commutative. 
Let $\Om^{\bul}$ (resp.~ $\Om'{}^{\bul}$) be a dga over $A$  
(resp.~ $A'$) 
such that each $\Om^q$ (resp.~ $\Om'{}^q$) $(q\in {\mab Z})$ is a $B$-module 
(resp.~ $B'$-module).  
Let $\varphi \col M\lo {\rm Ker}(\Om^1\lo \Om^2)$ 
(resp.~$\varphi' \col M\lo {\rm Ker}(\Om'{}^1\lo \Om'{}^2)$) be 
a morphism of $A$-modules (resp.~ $A'$-modules).  
Let $h\col \Om^{\bul}\lo \Om'{}^{\bul}$ be a morphism of complexes of $A$-modules 
fitting into the following commutative diagram 
\begin{equation*} 
\begin{CD} 
M@=M\\
@V{\varphi}VV @VV{\varphi'}V\\
\Om^1@>{h}>>\Om'{}^1. 
\end{CD} 
\end{equation*}
Let $E$ (resp.~$E'$) be a $B$-module (resp.~$B'$-module). 
Let $g\col E\lo E'$ be a morphism of $B$-modules. 
(We can consider $E'$ as a $B$-module.)
Let $(E\otimes_B\Om^{\bul}, d)$ (resp.~$(E'\otimes_{B'}\Om'{}^{\bul}, d')$) 
be a complex of $A$-modules (resp.~a complex of $A'$-modules).  
Assume that $f:=g\otimes h\col E\otimes_B\Om^{\bul}\lo E'\otimes_{B'}\Om'{}^{\bul}$ 
is a morphism of complexes of $A$-modules. 
Then we have the following morphism 
\begin{align*}  
f_{\langle M\rangle} \col 
(E\otimes_{B}\Om^{\bul})\langle M\rangle \lo (E'\otimes_{B'}\Om'^{\bul})\langle M\rangle. 
\end{align*}  
of complexes. Consider the mapping cone ${\rm MC}(f_{\langle M\rangle})$ of 
the morphism $f_{\langle M\rangle}$. 
Then 
\begin{align*} 
{\rm MC}(f_{\langle M\rangle})^q&= \{\Gam_A(M)\otimes_A(E\otimes_{B}\Om^{q+1}) \}
\oplus \{\Gam_A(M)\otimes_A(E'\otimes_{B'}\Om'^{q})\}\\
&=\Gam_A(M)\otimes_A(E\otimes_{B}\Om^{q+1}\oplus E'\otimes_{B'}\Om'^{q})\\
&=\Gam_A(M)\otimes_A{\rm MC}(f)^q. 
\end{align*} 
\par  
The following is the commutativity of the operation of the mapping cone and that of 
the PD-Hirsch extension. 

\begin{lemm}[{\rm {\bf \cite[(3.11)]{nhi}}}]\label{lemm:mcfgq}
The following diagram is commutative$:$
\begin{equation*} 
\begin{CD} 
{\rm MC}(f_{\langle M\rangle})^q
@=
\Gam_A(M)\otimes_A{\rm MC}(f)^q\\
@V{d_H,(\varphi,\varphi')}VV @VV{d_H,(-\varphi, \varphi')}V \\
{\rm MC}(f_{\langle M\rangle})^{q+1}
@=
\Gam_A(M)\otimes_A{\rm MC}(f)^{q+1}. 
\end{CD}
\tag{2.1.1}\label{cd:mqq}
\end{equation*} 
$($Note that the delicate sign appears before $\varphi)$. 
\end{lemm}

We have also obtained the following: 

\begin{prop}[{\rm {\bf \cite[(3.12)]{nhi}}}]\label{prop:acy}
Let the notations be as above. 
Assume that $M$ is a direct summand of a free $A$-module of countable rank. 
If $(\bigoplus_{i=1}^n(E_i\otimes_{B_i}\Om^{\bul}_i),d)$ is acyclic, 
then $(\bigoplus_{i=1}^n(E_i\otimes_{B_i}\Om^{\bul}_i))\langle M\rangle$ is acyclic. 
\end{prop}

\par
By using (\ref{prop:acy}),  
we obtain the following result which we need in this article:

\begin{coro}[{\rm {\bf \cite[(3.13)]{nhi}}}]\label{coro:hac}
Let $f \col E\otimes_B\Om^{\bul}\lo E'\otimes_{B'}\Om'^{\bul}$ be 
a quasi-isomorphism of $A$-modules. 
Assume that $M$ is a direct summand of a free $A$-module of countable rank. 
Then the morphism 
$f_{\langle M\rangle} \col (E\otimes_{B}\Om^{\bul})\langle M\rangle 
\lo (E'\otimes_{B'}\Om'^{\bul})\langle M\rangle$   
is a quasi-isomorphism of $A$-modules. 
\end{coro} 

\begin{rema}
In the later sections, we consider the Hirsch extension in characteristic $0$. 
Hence the notion of the PD-Hirsch extension and that of the Hirsch extension 
is the same. In the later sections we denote 
$\langle M\rangle$ by $[M]$. 
\end{rema}

\section{The key result}\label{sec:da}
In this section we give a key result (\ref{coro:fdc}) below in this article. 
\par 
Let $r$ be a positive integer. 
Let $S$ be a log analytic family over ${\mab C}$ of log points of virtual dimension $r$,  
that is, zariski locally on $S$, the log structure $M_S$ of $S$ is isomorphic to  
${\mab N}^r\oplus {\cal O}_S^*\owns ((n_1,\ldots,n_r),a)\lom 0^{\sum_{i=1}^rn_i}a\in  
{\cal O}_S$ (cf.~\cite[\S2]{nb}), where $0^n=0$ when $n\not =0$ and $0^0:=1$. 
Set $\ol{M}_S:=M_S/{\cal O}_S^*$. 
Let $g\col Y\lo S$ be a log smooth morphism of log analytic spaces over ${\mab C}$. 
\par
Locally on $S$,  there exists a family  
$\{t_i\}_{i=1}^r$ of local sections of $\ol{M}_S$ 
giving a local basis of $\ol{M}_S$. 
Let $M_i$ be the ideal sheaf of $M_S$ generated by a lift of $t_i$ in $M_S$. 
For all $1\leq \forall i\leq r$, the submonoid sheaf 
$\bigoplus_{j=1}^iM_j$ of $M_S$ and $\os{\circ}{S}$ 
defines a family of log points of virtual dimension $i$. 
Let $S_i:=S(M_1,\ldots,M_i)=(\os{\circ}{S},(\bigoplus_{j=1}^iM_j\lo {\cal O}_S)^a)$ 
be the resulting local log analytic space. Set $S_0:=\os{\circ}{S}$. 
Then we have the following sequence of families of log points 
of virtual dimensions locally: 
\begin{align*} 
S=S_r\lo S_{r-1}\lo S_{r-2}\lo \cdots \lo S_1\lo S_0=\os{\circ}{S}.
\tag{3.0.1}\label{ali:ffaoy}
\end{align*} 
Let $\wt{t}_i$ be a lift of $t_i$ in $M_S$. 
The one-form $d\log \wt{t}_i\in \Om^1_{S/\os{\circ}{S}}$ is 
independent of the choice of the lift $\wt{t}_i$. 
Hence we can denote $d\log \wt{t}_i$ by $d\log t_i$. 
Denote also by $d\log t_i\in \Om^1_{Y/\os{\circ}{S}}$ 
the image of $d\log t_i\in \Om^1_{S/\os{\circ}{S}}$ in 
$\Om^1_{Y/\os{\circ}{S}}$.  
Let ${\cal F}$ be a (not necessarily coherent) locally free ${\cal O}_Y$-module and let 
\begin{align*} 
\nabla \col {\cal F}\lo {\cal F}\otimes_{{\cal O}_Y}\Om^1_{Y/\os{\circ}{S}}
\tag{3.0.2}\label{ali:ffoy}
\end{align*} 
be an integrable connection. 
Then we have the complex 
${\cal F}\otimes_{{\cal O}_Y}\Om^{\bul}_{Y/\os{\circ}{S}}$ 
of $g^{-1}({\cal O}_S)$-modules.


\begin{lemm}\label{lemm:qseux}
$(1)$ The sheaf $\Om^i_{Y/\os{\circ}{S}}$ $(i\in {\mab N})$ is a locally free 
${\cal O}_Y$-module. 
\par 
$(2)$ Locally on $S$, the following sequence 
\begin{align*} 
0  \lo {\cal F}
\otimes_{{\cal O}_Y}
{\Om}^{\bul}_{Y/S_i}[-1] 
\os{{\rm id}_{\cal F}\otimes (d\log t_i\wedge)}{\lo} 
{\cal F}\otimes_{{\cal O}_Y}
{\Om}^{\bul}_{Y/S_{i-1}}  
\lo {\cal F}\otimes_{{\cal O}_Y}{\Om}^{\bul}_{Y/S_i} \lo 0 
\quad (1\leq i\leq r) 
\tag{3.1.1}\label{ali:agaxd}
\end{align*} 
is exact. 
$($Note that the morphism 
\begin{align*} 
{\rm id}_{\cal F}\otimes (d\log t_i\wedge) \col {\cal F}
\otimes_{{\cal O}_Y}
{\Om}^{\bul}_{Y/S_i}[-1] \lo 
{\cal F}\otimes_{{\cal O}_Y}
{\Om}^{\bul}_{Y/S_{i-1}}  
\end{align*} 
is indeed a morphism of complexes of $g^{-1}({\cal O}_S)$-modules.$)$ 
\end{lemm}
\begin{proof} 
(1), (2): First note that 
the connection (\ref{ali:ffoy}) induces the local integrable connection 
${\cal F}\lo {\cal F}\otimes_{{\cal O}_Y}{\Om}^1_{Y/S_i}$ for any $0\leq i\leq r$.   
Consequently we indeed have the local log de Rham complex 
${\cal F}\otimes_{{\cal O}_Y}{\Om}^{\bul}_{Y/S_i}$. 
Consider the following exact sequence (cf.~\cite[(3.12)]{klog1}): 
\begin{align*} 
f^*({\Om}^1_{S/\os{\circ}{S}}) \lo 
\Om^1_{Y/\os{\circ}{S}}  \lo \Om^1_{Y/S} \lo 0. 
\tag{3.1.2}\label{ali:aggaxd}
\end{align*} 
Because $Y/S$ is log smooth, 
the following sequence is exact and locally split (cf.~[loc.~cit.]): 
\begin{align*} 
0\lo f^*(\Om^1_{S/\os{\circ}{S}}) \lo 
\Om^1_{Y/\os{\circ}{S}}  \lo \Om^1_{Y/S} \lo 0. 
\tag{3.1.3}\label{ali:agfaxd}
\end{align*} 
Hence we obtain (1). 
\par 
Since $\{d\log t_i\}_{i=1}^r$ is a basis of $\Om^1_{S/\os{\circ}{S}}$, 
we see that $\{d\log t_i\}_{i=1}^r$ can be a part of a basis of 
$\Om^1_{Y/\os{\circ}{S}}$ by (\ref{ali:agfaxd}). 
Hence the following sequence 
\begin{align*} 
0  \lo {\cal O}_{Y} \os{d\log t_i\wedge}{\lo} \Om^1_{Y/S_{i-1}}  
\lo \Om^1_{Y/S_i} \lo 0 \quad (1\leq i\leq r)
\tag{3.1.4}\label{ali:agaddxd}
\end{align*} 
is a locally spit exact sequence.  
Now it is clear that the sequence (\ref{ali:agaxd}) is exact. 
\end{proof} 

\begin{rema}
In the exact sequence in \cite[(3.12)]{klog1}, $\om^1_{Y/S}$ 
(resp.~$\om^1_{X/S}$) should be replaced by 
$\om^1_{X/S}$ 
(resp.~$\om^1_{Y/S}$).
\end{rema}


\par 
Let $\{t_1,\ldots,t_r\}$ be a set of local generators of $\ol{M}_S$. 
Set 
$$U_S(t_1,\ldots,t_r):=\bigoplus_{i=1}^r{\cal O}_St_i$$ 
(a free ${\cal O}_S$-module with basis $t_1,\ldots,t_r$). 
In the following we denote $t_i$ in $U_S$ by $u_i$. 
Let $\Gam_{{\cal O}_S}(U_S):={\rm Sym}_{{\cal O}_S}(U_S)={\cal O}_S[\ol{M}_S]$ 
be the symmetric algebra of $U_S$ over ${\cal O}_S$.  
Consider the following Hirsch extension 
\begin{equation*}  
{\cal F}\otimes_{{\cal O}_Y}
\Om^{\bul}_{Y/\os{\circ}{S}}[U_S]
:= \Gam_{{\cal O}_S}(U_S)\otimes_{{\cal O}_S}
{\cal F}\otimes_{{\cal O}_Y}
\Om^{\bul}_{Y/\os{\circ}{S}} 
\tag{3.2.2}\label{eqn:blhff}
\end{equation*} 
of ${\cal F}\otimes_{{\cal O}_Y}
\Om^{\bul}_{Y/\os{\circ}{S}}$ 
by the morphism 
\begin{align*} 
d\log \col g^{-1}(U_S)\owns u_i\lom d\log t_i\in 
{\rm Ker}(\Om^1_{Y/\os{\circ}{S}}\lo \Om^2_{Y/\os{\circ}{S}})
\tag{3.2.3}\label{eqn:blff}
\end{align*}  
(\cite[\S3]{nhi}). 
We have omitted to write $g^{-1}$ for  
$g^{-1}(\Gam_{{\cal O}_S}(U_S))\otimes_{g^{-1}({\cal O}_S)}$ 
in (\ref{eqn:blhff}).
The differential morphism 
\begin{equation*} 
\nabla \col  {\cal F}\otimes_{{\cal O}_Y}
\Om^i_{Y/\os{\circ}{S}}[U_S]
\lo  {\cal F}\otimes_{{\cal O}_Y}
\Om^{i+1}_{Y/\os{\circ}{S}}[U_S]
\quad (i\in {\mab Z}_{\geq 0})
\end{equation*}   
is defined by the following: 
\begin{align*} 
&\nabla(m_1^{[e_1]}\cdots m_r^{[e_r]}\otimes \om)=
\sum_{j=1}^rm_1^{[e_1]}\cdots m_j^{[e_j-1]} \cdots m_r^{[e_r]}
d\log m_j\wedge \om +
m_1^{[e_1]}\cdots  m_r^{[e_r]}\otimes \nabla(\om)
\tag{3.2.4}\label{eqn:bdff}\\
&\quad (~\om \in {\cal F}
\otimes_{{\cal O}_Y}\Om^j_{Y/\os{\circ}{S}},m_1, \ldots m_r\in U_S, 
e_1,\ldots,e_r\in {\mab Z}_{\geq 1}, m_i^{[e_i]}=((e_i)!)^{-1}m^{e_i}_i~). 
\end{align*}
We can easily see that $\nabla^2=0$. 
\par 
The projection $\Gam_{{\cal O}_S}(U_S)
\lo \Gam_{{\cal O}_S,0}(U_S)={\cal O}_S$
induces the following natural morphism of complexes: 
\begin{equation*} 
{\cal F}\otimes_{{\cal O}_Y}
\Om^{\bul}_{Y/\os{\circ}{S}}[U_S]
\lo {\cal F}\otimes_{{\cal O}_Y}
\Om^{\bul}_{Y/S}.
\tag{3.2.5}\label{ali:uafacui} 
\end{equation*}

The following is a simpler analytic version of \cite[(4.5)]{nhi}:

\begin{defi}[{\rm {\bf cf.~\cite[(4.5)]{nhi}}}]\label{defi:wal}
(1) Take a local basis of $\Om^1_{Y/\os{\circ}{S}}$ containing 
$\{d\log t_1,\ldots,d\log t_r\}$. 
Let $(d\log t_i)^*\col \Om^1_{Y/\os{\circ}{S}}\lo {\cal O}_{Y}$ 
be the local morphism defined by the local dual basis of $d\log t_i$. 
We say that the connection $({\cal F},\nabla)$ has {\it no poles along} $S$ 
if the composite morphism 
${\cal F}\os{\nabla}{\lo} {\cal F}\otimes_{Y}\Om^1_{Y/\os{\circ}{S}}
\os{{\rm id}_{\cal F}\otimes(d\log t_i)^*}{\lo} {\cal F}$ vanishes for $1\leq \forall i\leq r$. 
\par 
(2) (cf.~\cite[p.~163]{kn}) 
If there exists locally a finite increasing filtration 
$\{{\cal F}_i\}_{i\in {\mab Z}}$ on ${\cal F}$ 
such that ${\rm gr}_i{\cal F}:={\cal F}_i/{\cal F}_{i-1}$ 
is a locally free ${\cal O}_Y$-module 
and $\nabla$ induces a connection 
${\cal F}_i\lo {\cal F}_i\otimes_{{\cal O}_Y}\Om^1_{Y/\os{\circ}{S}}$ 
whose induced connection 
${\rm gr}_i{\cal F}\lo ({\rm gr}_i{\cal F})
\otimes_{{\cal O}_Y}\Om^1_{Y/\os{\circ}{S}}$  has no poles along $M_S$, 
then we say that $({\cal F},\nabla)$ 
a {\it locally nilpotent integrable connection on $Y$ with respect to} $S$ 
\end{defi} 

It is obvious that 
the notion ``no poles along $M_S$'' 
is independent of the choice of the local sections $t_1,\ldots,t_r$.
It is also obvious that, if the connection (\ref{ali:ffoy}) factors through 
${\cal F}\otimes_{{\cal O}_Y}\Om^1_{\os{\circ}{Y}/\os{\circ}{S}}$, 
then $({\cal F}, {\nabla})$ has no poles along $S$. 
Especially $({\cal O}_Y^{\oplus n},d^{\oplus n}_Y)$ $(n\in {\mab Z}_{\geq 1})$ 
has no poles along $M_S$. 
Here $d_Y\col {\cal O}_Y\lo \Om^1_{Y/\os{\circ}{S}}$ is the usual derivative.

\par 
The following results (\ref{theo:r1}) and (\ref{theo:afp}) 
play the most important role in this article: 

\begin{theo}[{\rm {\bf cf.~\cite[(4.8)]{nhi}}}]\label{theo:r1}
Assume that $r=1$. Let $({\cal F},\nabla)$ be 
a locally nilpotent integrable connection on $Y$ with respect to $S$. 
Then the morphism {\rm (\ref{ali:uafacui})} is a quasi-isomorphism. 
\end{theo}
\begin{proof} 
Though the same proof as that of \cite[(4.8)]{nhi} works in the analytic case, 
we give the proof here for the completeness of this article.
\par 
This is a local problem. We may assume that $U_{S}={\cal O}_Su_1$. 
Denote $t_1$ and $u_1$ by $t$ and $u$, respectively. 
Assume that there exists a sub free ${\cal O}_Y$-module ${\cal F}_0$ of ${\cal F}$ 
such that ${\cal F}_0$ and ${\rm gr}_1{\cal F}:={\cal F}/{\cal F}_0$ 
have no poles along $S$ and ${\rm gr}_1{\cal F}$ is a free ${\cal O}_Y$-module. 
Let $\iota \col {\cal F}_0\otimes_{{\cal O}_Y}
\Om^{\bul}_{Y/\os{\circ}{S}} \os{\sus}{\lo} {\cal F}\otimes_{{\cal O}_Y}
\Om^{\bul}_{Y/\os{\circ}{S}}$ be the natural inclusion. 
For an injective morphism $f\col {\cal G}^{\bul}\lo {\cal H}^{\bul}$ of complexes of 
abelian sheaves on $\os{\circ}{Y}$, let ${\rm MC}(f)$ be the mapping cone of $f$. 
Then we have the natural quasi-morphism 
${\rm MC}(\iota)\os{\sim}{\lo} {\rm gr}_1{\cal F}\otimes_{{\cal O}_Y}
\Om^{\bul}_{Y/\os{\circ}{S}}$. 
By (\ref{coro:hac}) we have a quasi-isomorphism 
\begin{align*} 
{\rm MC}(\iota)[U_S]  \os{\sim}{\lo} {\rm gr}_1{\cal F}\otimes_{{\cal O}_Y}
\Om^{\bul}_{Y/\os{\circ}{S}}[U_S]. 
\end{align*} 
Hence, by (\ref{lemm:mcfgq}), we obtain the following isomorphism
\begin{align*} 
{\rm MC}(\iota_{[U_S]}\col {\cal F}_0\otimes_{{\cal O}_Y}
\Om^{\bul}_{Y/\os{\circ}{S}}[U_S]\lo 
{\cal F}\otimes_{{\cal O}_Y}\Om^{\bul}_{Y/\os{\circ}{S}}[U_S])  
\os{\sim}{\lo} {\rm gr}_1{\cal F}\otimes_{{\cal O}_Y}
\Om^{\bul}_{Y/\os{\circ}{S}}[U_S]. 
\end{align*} 
This means that the upper sequence of the following commutative diagram 
is exact: 
\begin{equation*}
\begin{CD}
0@>>> {\cal F}_0\otimes_{{\cal O}_Y}
\Om^{\bul}_{Y/\os{\circ}{S}}[U_S]@>>>
{\cal F}\otimes_{{\cal O}_Y}
\Om^{\bul}_{Y/\os{\circ}{S}}[U_S]@>>>
{\rm gr}_1{\cal F}\otimes_{{\cal O}_Y}
\Om^{\bul}_{Y/\os{\circ}{S}}[U_S]\lo 0\\
@. @VVV @VVV @VVV\\
0@>>> {\cal F}_0\otimes_{{\cal O}_Y}
\Om^{\bul}_{Y/S}@>>>
{\cal F}\otimes_{{\cal O}_Y}
\Om^{\bul}_{Y/S}@>>>
{\rm gr}_1{\cal F}\otimes_{{\cal O}_Y}
\Om^{\bul}_{Y/S}\lo 0. 
\end{CD}
\tag{3.4.1}\label{cd:uafecui} 
\end{equation*} 
Hence, if the morphism  
\begin{equation*} 
{\cal E}\otimes_{{\cal O}_Y}
\Om^{\bul}_{Y/\os{\circ}{S}}[U_S]
\lo {\cal E}\otimes_{{\cal O}_Y}
\Om^{\bul}_{Y/S}.
\end{equation*}  
is an isomorphism for ${\cal E}={\cal F}_0$ and ${\rm gr}_1{\cal F}$, 
then the morphism 
\begin{equation*} 
{\cal F}\otimes_{{\cal O}_Y}\Om^{\bul}_{Y/\os{\circ}{S}}[U_S]
\lo {\cal F}\otimes_{{\cal O}_Y}\Om^{\bul}_{Y/S}
\end{equation*}  
is an isomorphism. 
Consequently we may assume that ${\cal F}$ has no poles along $S$. 
By (\ref{ali:agaxd}) the following sequence is exact: 
\begin{align*} 
0  \lo {\cal F}\otimes_{{\cal O}_Y}\Om^{\bul}_{Y/S}[-1] \os{d\log t\wedge }{\lo} 
{\cal F}\otimes_{{\cal O}_Y}
{\Om}^{\bul}_{Y/\os{\circ}{S}}  
\lo {\cal F}\otimes_{{\cal O}_Y}\Om^{\bul}_{Y/S} \lo 0. 
\tag{3.4.2}\label{ali:agbxd}
\end{align*} 
Here we have denoted ${\rm id}_{\cal F}\otimes(d\log t\wedge)$ only by $d\log t\wedge$. 
We claim that the following sequence 
\begin{align*} 
0  \lo {\cal H}^{q-1}({\cal F}\otimes_{{\cal O}_Y}
{\Om}^{\bul}_{Y/S}) \os{d\log t\wedge }{\lo} 
{\cal H}^q({\cal F}\otimes_{{\cal O}_Y}
{\Om}^{\bul}_{Y/\os{\circ}{S}})  
\lo {\cal H}^q({\cal F}\otimes_{{\cal O}_Y}{\Om}^{\bul}_{Y/S}) \lo 0
\quad (q\in {\mab N}) 
\tag{3.4.3}\label{ali:agxhqd}
\end{align*} 
obtained by (\ref{ali:agbxd}) is exact. 
We have only to prove  
that the following inclusion holds:   
\begin{align*} 
\nabla({\cal F}\otimes_{{\cal O}_Y}
{\Om}^{q-1}_{Y/\os{\circ}{S}}) \cap 
(d\log t \wedge({\cal F}\otimes_{{\cal O}_Y}
{\Om}^{q-1}_{Y/\os{\circ}{S}})) \subset 
d\log t\wedge(\nabla ({\cal F}\otimes_{{\cal O}_Y}
{\Om}^{q-2}_{Y/\os{\circ}{S}})). 
\tag{3.4.4}\label{ali:dayn}
\end{align*} 
Set $\om_1:=d\log t$. 
Let us take a local basis $\{\om_1,\ldots,\om_m\}$ $(m\in {\mab N})$ of 
$\Om^1_{Y/\os{\circ}{S}}$ ((\ref{lemm:qseux})). 
Take a local section 
$$f:=\sum_{1\leq i_1<\cdots <i_{q-1}\leq m}
f_{i_1\cdots i_{q-1}}\om_{i_1}\wedge \cdots \wedge \om_{i_{q-1}}\quad 
(f_{i_1\cdots i_{q-1}}\in {\cal F})$$ 
of ${\cal F}\otimes_{{\cal O}_Y}
{\Om}^{q-1}_{Y/\os{\circ}{S}}$. 
Decompose $f$ as follows: 
$$f=\sum_{2\leq i_2<\cdots <i_{q-1}\leq m}f_{1 i_2\cdots i_{q-1}}\om_1
\wedge \om_{i_2} \wedge \cdots \wedge \om_{i_{q-1}}+
\sum_{2\leq i_1<\cdots <i_{q-1}\leq m}
f_{i_1\cdots i_{q-1}}\om_{i_1}\wedge \cdots \wedge \om_{i_{q-1}}.$$
Then 
\begin{align*} 
\nabla(f)&=\sum_{2\leq i_2< \cdots <i_{q-1}\leq m}\nabla(f_{1i_2\cdots i_{q-1}})
\wedge \om_1\wedge \om_{i_2}\wedge \cdots \wedge \om_{i_{q-1}}\tag{3.4.5}
\label{ali:tt}\\
&-
\sum_{2\leq i_2< \cdots <i_{q-1}\leq m}f_{1i_2\cdots i_{q-1}}\om_1
\wedge d(\om_{i_2}\wedge\cdots \wedge \om_{i_{q-1}}) \\
&+\sum_{2\leq i_1<\cdots <i_{q-1}\leq m}(\nabla
(f_{i_1\cdots i_{q-1}}\om_{i_1}\wedge \cdots \wedge \om_{i_{q-1}})).  
\end{align*}
Because ${\cal F}$ has no poles along $M_{S}$, 
the third term in (\ref{ali:tt}) can be expressed as follows:    
\begin{align*}
\sum_{2\leq i_1<\cdots <i_{q}\leq m}g_{i_1\cdots i_{q}}\om_{i_1}\wedge \cdots \wedge \om_{i_{q}}
\tag{3.4.6}\label{ali:ttvt}\\
\end{align*} 
for some $g_{i_1\cdots i_{q}}\in {\cal F}$.  
Assume that $\nabla(f)\in 
\om_1 \wedge({\cal F}\otimes_{{\cal O}_Y}{\Om}^{q-1}_{Y/\os{\circ}{S}})$. 
Then (\ref{ali:ttvt}) vanishes. 
Hence 
\begin{align*}
\nabla(f)&=-\om_1\wedge 
(\sum_{2\leq i_2<\cdots <i_{q-1}\leq m} \nabla(f_{1\cdots i_{q-1}}) 
\wedge \om_{i_2}\wedge \cdots \wedge \om_{i_{q-1}}\\
&+ 
\sum_{2\leq i_2<\cdots <i_{q-1}\leq m}f_{1\cdots i_{q-1}}
\wedge d(\om_{i_2}\wedge \cdots \wedge \om_{i_{q-1}}))\\
&=-\om_1\wedge 
\nabla(\sum_{2\leq i_2\cdots <i_{q-1}\leq m} f_{1\cdots i_{q-1}}\om_{i_2}\wedge
\cdots \wedge \om_{i_{q-1}}). 
\end{align*} 
Thus we have proved that the inclusion (\ref{ali:dayn}) holds. 
\par 
By (\ref{ali:agxhqd}) the following sequence is exact: 
\begin{align*} 
0  &\lo {\cal H}^0({\cal F}
\otimes_{{\cal O}_Y}
{\Om}^{\bul}_{Y/\os{\circ}{S}}) \os{d\log t\wedge }{\lo} 
{\cal H}^1({\cal F}
\otimes_{{\cal O}_Y}
{\Om}^{\bul}_{Y/\os{\circ}{S}})
 \os{d\log t\wedge }{\lo}  \cdots  \os{d\log t\wedge }{\lo}  
 {\cal H}^{q-1}({\cal F}\otimes_{{\cal O}_Y}
{\Om}^{\bul}_{Y/\os{\circ}{S}}) \tag{3.4.7}\label{ali:qf}\\
& \os{d\log t\wedge }{\lo} {\cal H}^q({\cal F}\otimes_{{\cal O}_Y}
{\Om}^{\bul}_{Y/\os{\circ}{S}}) 
\lo {\cal H}^q({\cal F}\otimes_{{\cal O}_Y}{\Om}^{\bul}_{Y/S}) \lo 0. 
\end{align*} 
\par 
\par 
Now we can complete the proof of this theorem by using the following lemma (\ref{lemm:lm}), 
which was proved in \cite{nhi}. 
\end{proof}

\begin{lemm}[\rm {\bf \cite[(4.9)]{nhi}}]\label{lemm:lm}
Let $A$ be a commutative ring with unit element. 
Let 
$$M^{\bul \bul}=(\prod_{i\in {\mab N},j\in {\mab Z}}M^{-i,j},
\prod_{i\in {\mab N},j\in {\mab Z}}d^{-i,j}, \prod_{i\in {\mab N},j\in {\mab Z}}d'^{-i,j})$$ 
be a double complex of $A$-modules, 
where $d^{-i,j}\col M^{-i,j}\lo M^{-i+1,j}$ and 
$d'{}^{-i,j}\col M^{-i,j}\lo M^{-i,j+1}$ are morphism of $A$-modules  
such that $d^{-i+1j}\circ d^{-i,j}=0$, $d'^{-i,j+1}\circ d'{}^{-i,j}=0$ and 
$d'{}^{-i+1,j}\circ d^{-i,j}+d^{-i,j+1}\circ d'{}^{-i,j}=0$. 
Let $N^{\bul}$ be a complex of $A$-modules. 
Let $\eps^j \col M^{0j}\lo N^j$ $(j\in {\mab Z})$ be a morphism of $A$-modules  
such that $\eps^{\bul}:=\{\eps^j\}_{j\in {\mab Z}}$ induces morphisms
$M^{0\bul}\lo N^{\bul}$ and $M^{-1,j}\os{d^{-1,j}}{\lo} M^{0j}\os{\eps^j}{\lo} N^j$ 
$(j\in {\mab Z})$ of complexes of $A$-modules.  
Assume that the following sequence 
\begin{align*}
\cdots \os{d^{-i-1,j}}\lo M^{-i,j}\os{d^{-i,j}}{\lo} 
M^{-i+1,j}\os{d^{-i+1,j}}{\lo} \cdots M^{0j}\os{\eps^j}{\lo} N^j\lo 0
\tag{3.5.1}\label{ali:nmij}
\end{align*} 
is exact and that 
the following natural sequence 
\begin{align*}
H^j(M^{-i,\bul},d'{}^{-i,\bul})\os{H^j(d^{-i,j})}{\lo} 
H^j(M^{-i+1,\bul},d'{}^{-i+1,\bul})\os{H^j(d^{-i+1,j})}{\lo} H^j(M^{-i+2,\bul},d'{}^{-i+2,\bul})
\tag{3.5.2}\label{ali:nhmij}
\end{align*}
is exact. 
Denote the sub-double complex 
$(\bigoplus_{i,j\in {\mab N}}M^{-i,j},\bigoplus_{i,j\in {\mab N}}d^{-i,j}, \bigoplus_{i,j\in {\mab N}}d'^{-i,j})_{i\in {\mab N},j\in {\mab Z}}$ 
of $M^{\bul \bul}$ 
simply by $\bigoplus_{i,j\in {\mab N}}M^{-i,j}$. 
Then 
the natural morphism $s(\bigoplus_{i,j\in {\mab N}}M^{-i,j})\lo N^{\bul}$ of complexes of $A$-modules 
induces the following injective morphism 
\begin{align*} 
H^q(s(\bigoplus_{i,j\in {\mab N}}M^{-i,j}))\os{\subset}{\lo} H^q(N^{\bul}) \quad (q\in {\mab Z}). 
\tag{3.5.3}\label{ali:nmmmij}
\end{align*} 
Furthermore, assume that the following natural morphism 
$${\rm Ker}(M^{0j}\lo M^{0,j+1})\lo {\rm Ker}(N^{0j}\lo N^{0,j+1})$$  
is surjective. 
Then the natural morphism $s(\bigoplus_{i,j\in {\mab N}}M^{-i,j})\lo N^{\bul}$ of complexes 
induces the following isomorphism 
\begin{align*} 
H^q(s(\bigoplus_{i,j\in {\mab N}}M^{-i,j}))\os{\sim}{\lo} H^q(N^{\bul}) \quad (q\in {\mab Z}). 
\tag{3.5.4}\label{ali:nmsmij}
\end{align*} 
\end{lemm}

The following is a key result in this article: 

\begin{theo}\label{theo:afp}
Let $({\cal F},\nabla)$ be 
a locally nilpotent integrable connection on $Y$ with respect to $S$.  
Then the morphism {\rm (\ref{ali:uafacui})} is a quasi-isomorphism. 
\end{theo} 
\begin{proof} 
This is a local problem. Hence we may assume that there exists 
the sequence (\ref{ali:ffaoy}). 
Consider the exact sequence (\ref{ali:agaxd}) for the case $i=r$. 
By the same proof as that of (\ref{theo:r1}), 
the following sequence 
\begin{align*} 
0 & \lo {\cal H}^{q-1}({\cal F}
\otimes_{{\cal O}_Y}
{\Om}^{\bul}_{Y/S}) \os{d\log t_r\wedge }{\lo} 
{\cal H}^q({\cal F}\otimes_{{\cal O}_Y}
{\Om}^{\bul}_{Y/S_{r-1}})  
&\lo {\cal H}^q({\cal F}\otimes_{{\cal O}_Y}{\Om}^{\bul}_{Y/S}) \lo 0 
\quad (q\in {\mab N}) \tag{3.6.1}\label{ali:agarqd}
\end{align*} 
obtained by (\ref{ali:agaxd}) is exact. 
Let 
\begin{align*} 
{\cal F}\otimes_{{\cal O}_Y}
\Om^{\bul}_{Y/S_{i-1}}[u_i,\ldots,u_j]:= \Gam_{{\cal O}_S}(\bigoplus_{k=i}^j
{\cal O}_Su_k)\otimes_{{\cal O}_S}
{\cal F}\otimes_{{\cal O}_Y}
\Om^{\bul}_{Y/S_{i-1}} \quad (1\leq i\leq j\leq r)
\end{align*}
be the Hirsch extension of ${\cal F}\otimes_{{\cal O}_Y}
\Om^{\bul}_{Y/S_{i-1}}$ by the following morphism 
$$\bigoplus_{k=i}^j{\cal O}_Su_k\owns u_k\lom d\log t_k\in 
\Om^1_{Y/S_{i-1}}.$$ 
By the same proof as that of (\ref{theo:r1}) again, we see that 
the following morphism
\begin{align*} 
{\cal F}\otimes_{{\cal O}_Y}
\Om^{\bul}_{Y/S_{r-1}}
[u_r]  \lo 
{\cal F}\otimes_{{\cal O}_Y}
\Om^{\bul}_{Y/S}
\tag{3.6.2}\label{ali:aganrqd}
\end{align*}
is a quasi-isomorphism. 
Analogously, we have the following quasi-isomorphism 
\begin{align*} 
{\cal F}\otimes_{{\cal O}_Y}
\Om^{\bul}_{Y/S_{r-2}}[u_{r-1}] \os{\sim}{\lo} 
{\cal F}\otimes_{{\cal O}_Y}
\Om^{\bul}_{Y/S_{r-1}}. 
\tag{3.6.3}\label{ali:aprqd}
\end{align*}
By (\ref{ali:aganrqd}), (\ref{ali:aprqd}) and (\ref{coro:hac}) 
we see that the following morphism 
\begin{align*} 
({\cal F}\otimes_{{\cal O}_Y}
\Om^{\bul}_{Y/S_{r-2}}[u_{r-1}])[u_r]
\lo {\cal F}\otimes_{{\cal O}_Y}
\Om^{\bul}_{Y/S} 
\end{align*}
is a quasi-isomorphism. 
Since the source of this quasi-isomorphism is equal to 
${\cal F}\otimes_{{\cal O}_Y}
\Om^{\bul}_{Y/S_{r-2}}[u_{r-1},u_r]$, 
we have the following quasi-isomorphism  
\begin{align*} 
{\cal F}\otimes_{{\cal O}_Y}
\Om^{\bul}_{Y/S_{r-2}}[u_{r-1},u_r]
\os{\sim}{\lo} {\cal F}\otimes_{{\cal O}_Y}
\Om^{\bul}_{Y/S}. 
\end{align*}
Continuing this process, we see that 
the morphism {\rm (\ref{ali:uafacui})} is a quasi-isomorphism. 
\end{proof}

Let $F$ be the filtration on 
${\cal F}\otimes_{{\cal O}_Y}\Om^{\bul}_{Y/S}$ defined by the following formula: 
\begin{align*} 
F^i({\cal F}\otimes_{{\cal O}_Y}\Om^{\bul}_{Y/S}):=
{\cal F}\otimes_{{\cal O}_Y}\Om^{\bul \geq i}_{Y/S}.
\end{align*} 
Let $\Gam_{{\cal O}_S,n}(U_S)$
be the degree $n$-part of $\Gam_{{\cal O}_S}(U_S)$.  
Let $F$ be the filtration on 
${\cal F}\otimes_{{\cal O}_Y}\Om^{\bul}_{Y/\os{\circ}{S}}[U_S]$ 
defined by the following formula: 
\begin{align*} 
F^i({\cal F}\otimes_{{\cal O}_Y}
\Om^{\bul}_{Y/\os{\circ}{S}}[U_S]):=
\bigoplus_{n\geq 0}\Gam_{{\cal O}_S,n}(U_S)\otimes_{{\cal O}_S} 
{\cal F}\otimes_{{\cal O}_Y}\Om^{\bul \geq i-n}_{Y/\os{\circ}{S}}.
\end{align*} 
It is clear that $F$ indeed gives the filtration on the complex 
${\cal F}\otimes_{{\cal O}_Y}
\Om^{\bul}_{Y/\os{\circ}{S}}[U_S]$ by (\ref{eqn:bdff}). 
Note that 
the filtration $F$ is separated, but {\it not} finite. 
The morphism (\ref{ali:uafacui}) induces the following filtered morphism 
\begin{align*} 
({\cal F}\otimes_{{\cal O}_Y}\Om^{\bul}_{Y/\os{\circ}{S}}[U_S],F)
\lo 
({\cal F}\otimes_{{\cal O}_Y}\Om^{\bul}_{Y/S},F).
\tag{3.6.4}\label{ali:aprqyd}
\end{align*} 

The following is a  generalization of \cite[(5.3)]{fut}: 

\begin{prop}[{\rm {\bf cf.~\cite[(5.3)]{fut}}}]\label{prop:anf}
The following morphism   
\begin{align*} 
{\rm gr}_F^i({\cal F}\otimes_{{\cal O}_{{Y}}}\Om^{\bul}_{{Y}/\os{\circ}{S}}[U_S])
\lo 
{\rm gr}_F^i({\cal F}\otimes_{{\cal O}_{{Y}}}\Om^{\bul}_{{Y}/S})
\tag{3.7.1}\label{ali:fggr}
\end{align*} 
induced by {\rm (\ref{ali:aprqyd})} is a quasi-isomorphism.  
\end{prop} 
\begin{proof} 
(The following proof is essentially the same as that of \cite[(5.3)]{fut}.)
First recall the definition of the Koszul complex. 
\par
For a ringed topos $({\cal T},{\cal A})$ and a morphism 
$\psi \col {\cal E}_1\lo {\cal E}_2$ 
of ${\cal A}$-modules and a nonnegative integer $n$, 
${\rm Kos}(\psi,n)^q:=
\Gam_{{\cal A},n-q}({\cal E}_1)\otimes_{\cal A}\bigwedge^{q}{\cal E}_2$ 
for $q\geq 0$ with the following boundary morphism
\begin{align*} 
{\rm Kos}(\psi,n)^q\owns e_1^{[i_1]}\cdots e_k^{[i_k]}\otimes f\lom 
\sum_{j=1}^ke_1^{[i_1]}\cdots 
e_1^{[i_j-1]}\cdots e_k^{[i_k]}\otimes \psi(e_j)\wedge f\in {\rm Kos}(\psi,n)^{q+1},
\end{align*}  
where $e_1,\ldots, e_k$ and $f$ are local sections of 
${\cal E}_1$ and ${\cal E}_2$, respectively, 
and $i_1,\ldots,i_k$ are positive integers such that $\sum_{j=1}^ki_j=n-q$. 
Here $\Gam_{\cal A}({\cal E}_1)=\bigoplus_{m=0}^{\inf}\Gam_{{\cal A},m}({\cal E}_1)$ 
is the graded PD-algebra over ${\cal A}$ generated by ${\cal E}_1$. 
\par
Consider the following injective morphism 
\begin{align*} 
\varphi \col {\cal O}_Y\otimes_{{\cal O}_S}U_S 
\owns b\otimes \sum_{i=1}^ra_iu_i \lom b\otimes \sum_{i=1}^ra_ig^*(d\log t_i)\in 
\Om^1_{{Y}/\os{\circ}{S}}.  
\end{align*} 
Note that we do not consider the following morphism 
\begin{align*} 
{\rm id}_{\cal F}\otimes_{{\cal O}_S}\varphi \col 
{\cal F}\otimes_{{\cal O}_S}U_S \lo {\cal F}
\otimes_{{\cal O}_{{Y}}}\Om^1_{{Y}/\os{\circ}{S}}. 
\end{align*} 
The following sequence 
\begin{align*} 
0  \lo{\cal O}_Y\otimes_{{\cal O}_S}U_S  \os{\varphi}{\lo} 
\Om^1_{{Y}/\os{\circ}{S}}  \lo \Om^1_{{Y}/S} \lo 0
\tag{3.7.2}\label{eqn:gsftxd}
\end{align*} 
is exact. 
Hence ${\rm Coker}(\varphi)$ is a flat ${\cal O}_Y$-module since
$Y/S$ is log smooth. 
\par 
We obtain the following equalities: 
\begin{align*} 
{\rm gr}_F^i({\cal F}\otimes_{{\cal O}_{{Y}}}\Om^{k}_{{Y}/\os{\circ}{S}}[U_S])&=
\bigoplus_{j\geq 0}\Gam_{{\cal O}_S,j}(U_S)\otimes_{{\cal O}_S}
{\cal F}\otimes_{{\cal O}_Y}{\rm gr}_F^{i-j}\Om^{k}_{{Y}/\os{\circ}{S}}  
\tag{3.7.3}\label{eqn:gsad}\\
&=\Gam_{{\cal O}_S,i-k}(U_S)
\otimes_{{\cal O}_S} 
{\cal F}\otimes_{{\cal O}_{{Y}}}
\bigwedge^{k}\Om^{1}_{{Y}/\os{\circ}{S}}   \\
&={\cal F}\otimes_{{\cal O}_{{Y}}}
\Gam_{{\cal O}_S,i-k}(U_S)\otimes_{{\cal O}_S}
\bigwedge^{k}\Om^{1}_{{Y}/\os{\circ}{S}}\\
&={\cal F}\otimes_{{\cal O}_{{Y}}}{\rm Kos}(\varphi,i)^k. 
\end{align*}
Since the morphism 
$${\rm gr}_F^i({\cal F}\otimes_{{\cal O}_{{Y}}}\Om^{\bul}_{{Y}/\os{\circ}{S}}[U_S])
\lo 
{\rm gr}_F^i({\cal F}\otimes_{{\cal O}_{{Y}}}\Om^{\bul}_{{Y}/\os{\circ}{S}}[U_S])$$
obtained by the second term on the right hand side of (\ref{eqn:bdff}) is $0$, we see that 
${\rm gr}_F^i({\cal F}\otimes_{{\cal O}_{{Y}}}\Om^{\bul}_{{Y}/\os{\circ}{S}}[U_S])$ 
is equal to 
${\cal F}\otimes_{{\cal O}_{{Y}}}{\rm Kos}(\varphi,i)$ by (\ref{eqn:gsad}). 
Because ${\rm Coker}(\varphi)$ is a flat ${\cal O}_Y$-module, 
we have the following isomorphism 
\begin{align*} 
\Gam_{{\cal O}_Y,i-q}({\rm Ker}(\varphi))\otimes_{{\cal O}_Y}
\bigwedge^q{\rm Coker}(\varphi)\os{\sim}{\lo} {\cal H}^q({\rm Kos}(\varphi,i))
\end{align*} 
by \cite[Proposition 4.4.1.6]{ic}. 
Hence 
\begin{align*} 
{\cal F}\otimes_{{\cal O}_Y}
\bigwedge^i{\rm Coker}(\varphi)\os{\sim}{\lo} {\cal F}\otimes_{{\cal O}_Y}
{\cal H}^i({\rm Kos}(\varphi,i))={\cal H}^i({\cal F}\otimes_{{\cal O}_Y}{\rm Kos}(\varphi,i))
\tag{3.7.4}\label{ali:kip}
\end{align*} 
and 
\begin{align*} 
{\cal H}^q({\cal F}\otimes_{{\cal O}_Y}
{\rm Kos}(\varphi,i))=
{\cal F}\otimes_{{\cal O}_Y}
{\cal H}^q({\rm Kos}(\varphi,i))=0\quad (q\not=i).  
\tag{3.7.5}\label{ali:kqp}
\end{align*} 
On the other hand, the target of the morphism (\ref{ali:fggr}) is equal to
${\cal F}\otimes_{{\cal O}_{{Y}}}\Om^{i}_{Y/S}
={\cal F}\otimes_{{\cal O}_{{\cal O}_Y}}\bigwedge^i{\rm Coker}(\varphi)$.  
Hence the equalities (\ref{ali:kip}) and (\ref{ali:kqp}) 
tell us that the morphism (\ref{ali:fggr}) is a quasi-isomorphism. 
We can complete the proof of (\ref{prop:anf}).  
\end{proof} 

\begin{rema}
Because the filtration $F$ on 
${\cal F}\otimes_{{\cal O}_{{Y}}}\Om^{\bul}_{{Y}/\os{\circ}{S}}[U_S]$ 
is not finite, (\ref{prop:anf}) does not imply (\ref{theo:afp}). 
\end{rema}

\begin{coro}\label{coro:fdc}
The filtered morphism {\rm (\ref{ali:aprqyd})} is a filtered quasi-isomorphism 
in the sense of {\rm \cite{nh2}}. 
\end{coro}
\begin{proof} 
This follows from (\ref{theo:afp}) and (\ref{prop:anf}). 
\end{proof} 

\section{Main Result}\label{sec:mr}
In this section we give the main result in this article. 
\par 
Let the notations be as in the previous section. 
Let $f\col X\lo S$ be a log smooth morphism of analytic spaces.  
\par 
In the following we assume that there exists a 
set $\{\os{\circ}{X}_{\lam}\}_{\lam \in \Lam}$ of smooth closed analytic spaces  
of $\os{\circ}{X}$ such that $\bigcup_{\lam\in \Lam}\os{\circ}{X}_{\lam}=\os{\circ}{X}$ 
and such that, for any point $x\in \os{\circ}{X}$,  
there exists no empty open subset $\os{\circ}{V}$ of $\os{\circ}{X}$ containing $x$ 
such that $\os{\circ}{X}_{\lam}\vert_{\os{\circ}{V}}\not=
\os{\circ}{X}_{\lam'}\vert_{\os{\circ}{V}}$ for any $\lam\not=\lam'$ 
except in the case where both hand sides 
are empty sets.

\par 
Set $\Del^r_{\os{\circ}{S}}:=\Del^r\times \os{\circ}{S}$. 
Then the standard coordinates of $\Del^r$ define an fs log structure 
on $\Del^r_{\os{\circ}{S}}$. Let $\Del^r_{S}$ be the resulting log analytic space. 
Locally on $S$, we have the following natural exact closed immersion
\begin{align*} 
S\os{\sus}{\lo} \Del^r_{S}
\end{align*} 
of log analytic spaces. 
This immersion fits into the following commutative diagram
\begin{equation*} 
\begin{CD} 
S@>{\sus}>> \Del^r_{S}\\
@VVV @VVV \\
\os{\circ}{S}@=\os{\circ}{S},
\end{CD} 
\end{equation*} 
where the vertical morphisms are natural morphisms. 
\par
By replacing ${\mab A}^r_{S}$ in \cite[(1.1)]{nb} by $\Del^r_{S}$, 
we have a well-defined log analytic space $\ol{S}$ 
with an exact closed immersion $S\os{\sus}{\lo} \ol{S}$ over $\os{\circ}{S}$. 
Locally on $S$, $\ol{S}$ is isomorphic to a polydisk $\Del^r_S$. 
We also assume that, locally on $X$, the morphism $\os{\circ}{f}$ fits into 
the following cartesian diagram 
\begin{equation*} 
\begin{CD} 
\os{\circ}{X}@>{\sus}>> \Del^n_{\os{\circ}{S}}\\
@V{\os{\circ}{f}}VV @VVV \\
\os{\circ}{S}@>{\sus}>> \os{\circ}{\ol{S}}
\end{CD} 
\end{equation*} 
such that $\os{\circ}{X}_{\lam}$ is empty or defined  
by equations $x_{i_1}=\cdots =x_{i_k}=0$ 
for some $1\leq i_1<\cdots <i_k\leq n$, where $x_1,\ldots, x_n$ are 
the standard coordinates of $\Del^n_{\os{\circ}{S}}$ in this cartesian diagram.  
We also assume that $x_{i_1},\ldots, x_{i_k}$ of the images of local sections of 
the log structure of $X$ by the structural morphism. 

\begin{rema}[{\bf {cf.~\cite[(4.3)]{fut}}}]
In the following we do not use the following fact. 
\par 
Let $Y/S$ be as in the previous section. 
Locally on $Y$, there exists a log smooth lift 
$\ol{Y}$ into a log analytic space over $\ol{S}$ of $Y$ 
fitting into the following cartesian diagram 
\begin{equation*} 
\begin{CD} 
Y@>{\subset}>>\ol{Y}\\
@VVV @VVV\\
S@>{\subset}>>\ol{S}\\
\end{CD}
\tag{4.1.1}\label{eqn:lll}
\end{equation*}  
(cf.~\cite[(1.1.8)]{nb}). 
Hence the condition of the existence of $(U,M_U)$ in \cite[(4.3)]{fut} 
is automatically  satisfied (if $Y/S$ is log smooth). In the case above, $M_U$ is fine. 
\end{rema}

\par 
For a subset $\ul{\lam}=\{\lam_1,\ldots,\lam_m\}$ of $\Lam$, 
set $\os{\circ}{X}_{\ul{\lam}}
:=\os{\circ}{X}_{\lam_1}\cap \cdots \cap \os{\circ}{X}_{\lam_m}$ 
and let $\os{\circ}{\iota}_{\ul{\lam}}\col \os{\circ}{X}_{\ul{\lam}}\lo 
\os{\circ}{X}$ be the natural closed immersion. 
Endow $\os{\circ}{X}_{\ul{\lam}}$ with the inverse image of the log structure of $X$ 
by the morphism $\os{\circ}{\iota}_{\ul{\lam}}$. 
Let $X_{\ul{\lam}}$ be the resulting log analytic space.  
Then the structural morphism $X_{\ul{\lam}}\lo S$ is ideally log smooth. 
Here we say that a morphism $Y\lo S$ of fine log analytic spaces is 
{\it ideally log smooth} if, locally on $Y$, there exists an exact closed immersion 
$Y\os{\sus}{\lo} Z$ over $S$ into a log smooth analytic space such that the ideal of definition of this exact closed immersion is 
defined by the image of an ideal of the log structure of $Z$. 
Let $\iota_{\ul{\lam}}\col X_{\ul{\lam}}\lo X$ be the resulting exact closed immersion. 
Let $\Om^{\bul}_{X_{\ul{\lam}}/\os{\circ}{S}}$ and 
$\Om^{\bul}_{X_{\ul{\lam}}/S}$ be the log de Rham complex of  
$X_{\ul{\lam}}/\os{\circ}{S}$ and $X_{\ul{\lam}}/S$, respectively. 
It is convenient to define 
$X_{\emptyset}=X$ and $\iota_{\emptyset}={\rm id}_X$. 
In \cite{fut} Fujisawa has essentially proved the following two results: 

\begin{prop}[{\rm {\bf cf.~\cite[(4.13)]{fut}}}]\label{prop:lm}
For a subset $\ul{\lam}$ of $\Lam$,  
$\Om^i_{X_{\ul{\lam}}/\os{\circ}{S}}=
\iota_{\ul{\lam}}^*(\Om^i_{X/\os{\circ}{S}})$ $(i\in {\mab N})$ 
and 
$\Om^i_{X_{\ul{\lam}}/S}=
\iota_{\ul{\lam}}^*(\Om^i_{X/S})$ $(i\in {\mab N})$.  
Consequently $\Om^i_{X_{\ul{\lam}}/S}$ and $\Om^i_{X_{\ul{\lam}}/\os{\circ}{S}}$ 
are locally free ${\cal O}_{X_{\ul{\lam}}}$-modules. 
\end{prop}
\begin{proof} 
Let ${\cal I}_{\ul{\lam}}$ be the defining sideal sheaf of $X_{\ul{\lam}}$ of $X$. 
Set $T:=\os{\circ}{S}$ or $S$. 
Then we have the following exact sequence 
\begin{align*} 
{\cal I}_{\ul{\lam}}/{\cal I}_{\ul{\lam}}^2\lo 
{\cal O}_{X_{\ul{\lam}}}\otimes_{{\cal O}_X}\Om^1_{X/T}\lo 
\Om^1_{X_{\ul{\lam}}/T}\lo 0
\end{align*}
by \cite[(3.6) (1)]{kn}. 
By the assumption, the left morphism is zero as in [loc.~cit.~(3.6) (2)]. 
\end{proof}

\begin{prop}[{\rm {\bf cf.~\cite[(4.16)]{fut}}}]\label{prop:lkm}
Set $X^{(n)}:=\coprod_{\sharp \ul{\lam}=n+1}X_{\ul{\lam}}$ and 
let $a^{(n)}\col X^{(n)}\lo X$ be the natural morphism. 
Fix a total order on $\Lam$. For a subset 
$\ul{\lam}=\{\lam_0,\ldots,\lam_n\}$ $(\lam_0<\cdots <\lam_n)$ of $\Lam$, 
set $\ul{\lam}_j:=\{\lam_0,\ldots,\lam_{j-1},\lam_{j+1},\ldots,\lam_n\}$. 
Let 
$\iota^{\ul{\lam}_j}_{\ul{\lam}}\col X_{\ul{\lam}}\os{\sus}{\lo} X_{\ul{\lam}_j}$ 
be the natural inclusion. 
Let 
\begin{align*} 
\rho^{n-1} \col a^{(n-1)}_*({\cal O}_{X^{(n-1)}})\lo a^{(n)}_*({\cal O}_{X^{(n)}})
\tag{4.3.1}\label{ali:oomxx}
\end{align*} 
be the morphism 
$$\sum_{\sharp \ul{\lam}=n+1}(-1)^j\iota^{\ul{\lam}_j*}_{\ul{\lam}}
\col a^{(n-1)}_*(\bigoplus_{\sharp \ul{\mu}=n}
{\cal O}_{X_{\ul{\mu}}})\lo a^{(n)}_*(\bigoplus_{\sharp \ul{\lam}=n+1}
{\cal O}_{X_{\ul{\lam}}}).$$ 
Then the following sequence 
\begin{align*} 
0\lo {\cal O}_X\lo a^{(0)}_*({\cal O}_{X^{(0)}})\os{\rho^{0}}{\lo} 
a^{(1)}_*({\cal O}_{X^{(1)}})\os{\rho^{1}}{\lo}  a^{(2)}_*({\cal O}_{X^{(2)}})
\os{\rho^{2}}{\lo} \cdots
\tag{4.3.2}\label{ali:ooxx}
\end{align*} 
is exact. 
\end{prop}
\begin{proof}  
The exactness of (\ref{ali:ooxx}) is a local question. 
We have only to take the tensorization ${\cal O}_S\otimes_{\mab C}$ 
in the proof of \cite[(4.16)]{fut}. 
\end{proof} 

\begin{rema}
If  $\os{\circ}{X}$ is a locally product of SNC analytic spaces, 
we can prove (\ref{prop:lkm}) in a very quick way. 
Indeed this is well-known 
in the case where $\os{\circ}{X}$ is an SNC(=simple normal crossing) 
analytic space (cf.~\cite[p.~115]{st2}). 
In the general case,  because 
$\os{\circ}{X}$ is a locally product of SNC analytic spaces, 
this follows from the case of SNC analytic spaces. 
\end{rema}

\par 
In the rest of this article, we fix a total order on $\Lam$. 
Let ${\cal E}$ be a (not necessarily coherent) locally free ${\cal O}_X$-module and let 
\begin{align*} 
\nabla \col {\cal E}\lo {\cal E}\otimes_{{\cal O}_X}\Om^1_{X/\os{\circ}{S}}
\tag{4.4.1}\label{ali:ffioy}
\end{align*} 
be an integrable connection. This connection induces
the following integrable connection 
\begin{align*} 
{\cal E}\lo {\cal E}\otimes_{{\cal O}_X}\Om^1_{X/S}
\tag{4.4.2}\label{ali:ffloy}
\end{align*} 
and we have the log de Rham complex 
${\cal E}\otimes_{{\cal O}_X}\Om^{\bul}_{X/S}$. 
We also have the log de Rham complex  
${\cal E}\otimes_{{\cal O}_X}a^{(n)}_*(\Om^{\bul}_{X^{(n)}/S})$ 
for each $n\in {\mab N}$. 
Let $T$ be $\os{\circ}{S}$ or $S$.  
Then we have a morphism 
$$\rho^{n-1,j} \col {\cal E}\otimes_{{\cal O}_X}a^{(n-1)}_*(\Om^j_{X^{(n-1)}/T})
\lo {\cal E}\otimes_{{\cal O}_X}a^{(n)}_*(\Om^j_{X^{(n)}/T})$$  
as in (\ref{ali:oomxx}). 

\begin{coro}\label{coro:qq}
Let 
${\cal E}\otimes_{{\cal O}_X}a^{(\star)}_*(\Om^{\bul}_{X^{(\star)}/T})$ 
be the double complex whose $(i,j)$-component $K^{ij}$ is 
${\cal E}\otimes_{{\cal O}_X}a^{(i)}_*(\Om^j_{X^{(i)}/T})$ 
and 
whose horizontal boundary morphism 
is $\rho^{ij}\col K^{ij}\lo K^{i+1,j}$ and whose vertical morphism 
is $(-1)^i\nabla$. 
Let 
$s({\cal E}\otimes_{{\cal O}_X}a^{(\star)}_*(\Om^{\bul}_{X^{(\star)}/T}))$ 
be the single complex of 
${\cal E}\otimes_{{\cal O}_X}a^{(\star)}_*(\Om^{\bul}_{X^{(\star)}/T})$. 
Endow 
$s({\cal E}\otimes_{{\cal O}_X}a^{(\star)}_*(\Om^{\bul}_{X^{(\star)}/T}))$ 
with the Hodge filtration $F$ defined by the following formula: 
\begin{equation*} 
F^is({\cal E}\otimes_{{\cal O}_X}a^{(\star)}_*(\Om^{\bul}_{X^{(\star)}/T}))
:=s({\cal E}\otimes_{{\cal O}_X}a^{(\star)}_*(\Om^{\bul \geq i}_{X^{(\star)}/T}))
\quad (i\in {\mab N}). 
\tag{4.5.1}\label{ali:fise}
\end{equation*} 
Then the natural filtered morphism 
\begin{align*} 
({\cal E}\otimes_{{\cal O}_X}\Om^{\bul}_{X/T},F)
\lo 
(s({\cal E}\otimes_{{\cal O}_X}a^{(\star)}_*(\Om^{\bul}_{X^{(\star)}/T})),F)
\tag{4.5.2}\label{ali:fisexe}
\end{align*} 
is a filtered quasi-isomorphism. 
\end{coro}
\begin{proof}
This immediately follows from (\ref{prop:lm}) and (\ref{prop:lkm}). 
\end{proof}

Consider the Hirsch extension 
${\cal E}\otimes_{{\cal O}_X}\Om^{\bul}_{X/\os{\circ}{S}}[U_S]$ 
of ${\cal E}\otimes_{{\cal O}_X}\Om^{\bul}_{X/\os{\circ}{S}}$ 
as in the previous section. 
We also have the Hirsch extension 
${\cal E}\otimes_{{\cal O}_X}a_{\ul{\lam}*}
(\Om^{\bul}_{X_{\ul{\lam}}/\os{\circ}{S}})[U_S]$ 
of ${\cal E}\otimes_{{\cal O}_X}
a^{(n)}_*(\Om^{\bul}_{X^{(n)}/\os{\circ}{S}})[U_S]$ 
as in the previous section. 
We obtain the double complex 
${\cal E}\otimes_{{\cal O}_X}a^{(\star)}_*(\Om^{\bul}_{X^{(\star)}/\os{\circ}{S}})
[U_S]$
and the single complex 
$s({\cal E}\otimes_{{\cal O}_X}a^{(\star)}_*(\Om^{\bul}_{X^{(\star)}/\os{\circ}{S}})
[U_S])$ of it.

Set 
\begin{align*} 
H(X/S,{\cal E}):= 
s({\cal E}\otimes_{{\cal O}_X}a^{(\star)}_*(\Om^{\bul}_{X^{(\star)}/\os{\circ}{S}})
[U_S]). 
\end{align*}
Endow $H(X/S,{\cal E})$ with the Hodge filtration defined by 
the following formula 
\begin{align*} 
F^iH(X/S,{\cal E}):=
s(\bigoplus_{n\geq 0}
\Gam_{{\cal O}_S,n}(U_S)\otimes_{{\cal O}_S}{\cal E}
\otimes_{{\cal O}_X}a^{(\star)}_*
(\Om^{\bul \geq i-n}_{X^{(\star)}/\os{\circ}{S}})). 
\end{align*} 

\begin{coro}\label{coro:emm}
The following filtered morphism 
\begin{align*} 
({\cal E}\otimes_{{\cal O}_X}\Om^{\bul}_{X/\os{\circ}{S}}[U_S],F)
\lo 
(H(X/S,{\cal E}),F)
\tag{4.6.1}\label{ali:emu}
\end{align*} 
is a filtered quasi-isomorphism. 
\end{coro}
\begin{proof} 
By (\ref{coro:qq}) and (\ref{coro:hac}) we see that 
the following morphism 
\begin{align*} 
{\cal E}\otimes_{{\cal O}_X}\Om^{\bul}_{X/\os{\circ}{S}}[U_S]
\lo H(X/S,{\cal E})
\end{align*} 
is a quais-isomorphism. 
We have the following equalities: 
\begin{align*} 
{\rm gr}_F^i({\cal E}\otimes_{{\cal O}_X}\Om^{\bul}_{X/\os{\circ}{S}}[U_S])&:=
(\cdots \lo {\rm gr}_F^is(\bigoplus_{n\geq 0}\Gam_{{\cal O}_S,n}(U_S)
\otimes_{{\cal O}_S}{\cal E}
\otimes_{{\cal O}_X} 
(\Om^{k}_{X/\os{\circ}{S}}))\lo \cdots)\\
&=
(\cdots \lo s(\Gam_{{\cal O}_S,i-k}(U_S)\otimes_{{\cal O}_S}{\cal E}
\otimes_{{\cal O}_X}\Om^{k}_{X/\os{\circ}{S}}))\lo \cdots)
\end{align*} 
and 
\begin{align*} 
{\rm gr}_F^iH(X/S,{\cal E})&:=
({\rm gr}_F^is(\bigoplus_{n\geq 0}
{\rm Sym}^n_{{\cal O}_S}(U_S)\otimes_{{\cal O}_S}{\cal E}
\otimes_{{\cal O}_X}a^{(\star)}_*
(\Om^{k}_{X^{(\star)}/\os{\circ}{S}}))\lo \cdots)\\
&=
(\cdots \lo s(\Gam_{{\cal O}_S,i-k}(U_S)\otimes_{{\cal O}_S}{\cal E}
\otimes_{{\cal O}_X}a^{(\star)}_*
(\Om^{k}_{X^{(\star)}/\os{\circ}{S}}))\lo \cdots). 
\end{align*} 
Hence the natural morphism 
\begin{align*}
{\rm gr}_F^i({\cal E}\otimes_{{\cal O}_X}\Om^{\bul}_{X/\os{\circ}{S}}[U_S])
\lo 
{\rm gr}_F^iH(X/S,{\cal E})
\end{align*} 
is an isomorphism by the proof of (\ref{coro:qq}). 
\end{proof}

The following is the main result of this article: 

\begin{theo}\label{theo:mra}
There exists the following filtered isomorphism 
\begin{align*}
(H(X/S,{\cal E}),F)
\os{\sim}{\lo} ({\cal E}\otimes_{{\cal O}_X}\Om^{\bul}_{X/S},F) 
\tag{4.7.1}\label{ali:ksx}
\end{align*} 
in the derived category   
${\rm D}^+{\rm F}(f^{-1}({\cal O}_S))$ of bounded 
below filtered complexes of $f^{-1}({\cal O}_S)$-modules on $X$. 

\end{theo}
\begin{proof}
By (\ref{coro:fdc}) and (\ref{coro:emm}) 
we have the following filtered quasi-isomorphisms:  
\begin{equation*}
({\cal E}\otimes_{{\cal O}_X}\Om^{\bul}_{X/S},F)
\os{\sim}{\longleftarrow} 
({\cal E}\otimes_{{\cal O}_X}\Om^{\bul}_{X/\os{\circ}{S}}[U_S],F) 
\os{\sim}{\lo} (H(X/S,{\cal E}),F). 
\tag{4.7.2}\label{ali:ooas}
\end{equation*} 
Hence we have the filtered isomorphism (\ref{ali:ksx}). 
\end{proof}

\begin{rema}\label{rema:cddm}
By (\ref{coro:fdc}), (\ref{coro:qq}) and (\ref{coro:emm}) 
we have the following commutative diagram 
\begin{equation*}
\begin{CD}
({\cal E}\otimes_{{\cal O}_X}\Om^{\bul}_{X/\os{\circ}{S}}[U_S],F) 
@>{\sim}>>(H(X/S,{\cal E}),F)\\
@V{\simeq}VV @VVV \\
({\cal E}\otimes_{{\cal O}_X}\Om^{\bul}_{X/S},F)
@>{\sim}>> 
(s({\cal E}\otimes_{{\cal O}_X}a^{(\star)}_*(\Om^{\bul}_{X^{(\star)}/S})),F). 
\end{CD}
\tag{4.8.1}\label{ali:oonas}
\end{equation*} 
Hence the right vertical morphism in (\ref{ali:oonas}) is also a quasi-isomorphism. 
\end{rema}

The following is a generalization of (\ref{coro:nk}): 

\begin{coro}\label{coro:rfks}
There exists a filtered isomorphism 
\begin{align*} 
Rf_*((H(X/S,{\cal E}),F))\os{\sim}{\lo} 
Rf_*(({\cal E}\otimes_{{\cal O}_X}\Om^{\bul}_{X/S},F)). 
\tag{4.9.1}\label{ali:kxef}
\end{align*} 
\end{coro}

Lastly we consider the contravariant functoriality of the isomorphism (\ref{ali:ksx})
whose formulation is not obvious.
\par
Let $S'$ be an analytic family over ${\mab C}$ of log points of virtual dimension $r'$ 
and let $v\col S\lo S'$ be a morphism of log analytic spaces. 
This morphism induces a morphism $v^*\col \ol{M}_{S'}\lo v_*(\ol{M}_{S})$.  
Locally on $S$, this morphism  is equal to a morphism 
$v^*\col {\mab N}^{r'}\lo {\mab N}^r$. 
Set 
$e_i=(0, \ldots,0,\os{i}{1},0,\ldots, 0)\in {\mab N}^s$ for $s=r$ or $r'$  
$(1\leq i\leq r)$. 
Let $A=(a_{ij})_{1\leq i\leq r, 1\leq j\leq r'} \in M_{rr'}({\mab N})$ 
be the representing matrix of $v^*$: 
$$v^*(e_1,\ldots,e_{r'})=(e_1,\ldots,e_r)A.$$ 
Let $t_1,\ldots,t_r$ and $t'_1,\ldots, t'_{r'}$ be local sections of 
$M_S$ and $M_{S'}$ 
whose images in $\ol{M}_S\os{\sim}{\lo} {\mab N}^r$ and 
$\ol{M}{}'_{S'}\os{\sim}{\lo}{\mab N}^{r'}$ are 
$e_1,\ldots,e_r$ and $e_1,\ldots,e_{r'}$, respectively. 
Let $\ol{t}_i$ and $\ol{t}{}'_i$ be the images of $t_i$ and $t'_i$ in 
$\ol{M}_S$ and $\ol{M}_{S'}$, respectively. 
Then there exists a local section $b_j\in {\cal O}_{S}^*$ 
$(1\leq j\leq r')$ such that 
\begin{align*} 
v^*(t'_j)=b_jt^{a_{1j}}_1\cdots t^{a_{rj}}_r. 
\tag{4.9.2}\label{ali:vta} 
\end{align*}
Let $u_1,\ldots, u_r$ and $u'_1,\ldots,u'_{r'}$ be 
the corresponding local sections 
to $\ol{t}_1,\ldots, \ol{t}_r$ and $\ol{t}{}'_1,\ldots,\ol{t}{}'_{r'}$ of 
$U_S$ and $U_{S'}$, respectively. Then we define 
an ${\cal O}_{S'}$-linear morphism 
$v^*\col \Gam_{{\cal O}_{S'}}(U_{S'})\lo v_*(\Gam_{{\cal O}_S}(U_S))$
by the following formula: 
\begin{align*} 
v^*(u'_j)=a_{1j}u_1+\cdots+a_{rj}u_r.
\tag{4.9.3}\label{ali:vua} 
\end{align*}
and by ${\rm Sym}(v^*)$ for the $v^*$ in (\ref{ali:vua}). 
Since ${\rm Aut}({\mab N}^s)={\mathfrak S}_s$ (\cite[p.~47]{nh3}), 
the morphism $v^*\col \Gam_{{\cal O}_{S'}}(U_{S'})\lo v_*(\Gam_{{\cal O}_S}(U_S))$
is independent of the choice of the local isomorphisms 
$\ol{M}_{S'}\simeq {\mab N}^{r'}$ and $\ol{M}_{S}\simeq {\mab N}^r$. 
In conclusion, 
we obtain the following well-defined morphism 
\begin{align*} 
v^*\col \Gam_{{\cal O}_{S'}}(U_{S'})\lo 
v_*(\Gam_{{\cal O}_S}(U_S))
\tag{4.9.4}\label{ali:kxvef}
\end{align*} 
of sheaves of commutative rings of unit elements on $S'$.  
This morphism satisfies the usual transitive relation 
\begin{align*} 
{\textrm ``}(v\circ v')^*=
v'{}^*v^*{\textrm '}{\textrm '}.
\tag{4.9.5}\label{ali:kxbvef}
\end{align*} 

\begin{rema}
Note that the morphism (\ref{ali:kxvef}) is different from 
the natural morphism 
\begin{align*} 
{\cal O}_{S'}[\ol{M}_{S'}]\lo v_*{\cal O}_S[\ol{M}_S]
\tag{4.10.1}\label{ali:kxbivef}
\end{align*} 
induced by $v^*\col M_{S'}\lo v_*(M_S)$. 
It is clear that $U_S$ is more important than $\ol{M}_S$ in this article. 
\end{rema}

\par 
Let $f'\col X'\lo S'$ be an analogous morphism of log analytic spaces to 
$f\col X\lo S$. 
Let $X'_{\lam'}$'s be analogous log analytic spaces to $X_{\lam}$'s for $X'$.
Assume that we are given a commutative diagram 
\begin{equation*} 
\begin{CD}
X@>{g}>> X'\\
@V{f}VV @VV{f'}V \\
S@>{v}>> S'
\end{CD}
\tag{4.10.2}\label{ali:ooxas}
\end{equation*}
such that, for any smooth component 
$\os{\circ}{X}_{\lam}$ of $\os{\circ}{X}$ over $\os{\circ}{S}$, 
there exists a unique smooth component 
$\os{\circ}{X}{}'_{\lam'}$ of 
$\os{\circ}{X}{}'$ over $\os{\circ}{S}{}'$ such that $g$ 
induces a morphism 
$\os{\circ}{X}_{\lam} \lo \os{\circ}{X}{}'_{\lam'}$. 

\begin{prop}[{\bf Contravariant functoriality}]\label{prop:rks}
Assume that we are given a commutative diagram 
\begin{equation*} 
\begin{CD}
{\cal E}'@>{\nabla'}>> {\cal E}'\otimes_{{\cal O}_{X'}}\Om^1_{X'/\os{\circ}{S}{}'}\\
@VVV @VVV \\
g_*({\cal E}) @>{g_*(\nabla)}>> g_*({\cal E}\otimes_{{\cal O}_X}\Om^1_{X/\os{\circ}{S}}), 
\end{CD}
\tag{4.11.1}\label{cd:kgsx}
\end{equation*}
where $({\cal E}',{\nabla'})$ is an analogous connection to 
$({\cal E},{\nabla})$ on $\os{\circ}{X}{}'$. 
Then the filtered morphism {\rm (\ref{ali:ksx})} is 
contravariantly functorial for 
the commutative diagrams {\rm (\ref{ali:ooxas})} and {\rm (\ref{cd:kgsx})}. 
This contravariance satisfies the usual transitive relation ``$(g\circ h)^*
=h^*g^*$''. 
\end{prop} 
\begin{proof} 
By the assumption and the existence of the morphism (\ref{ali:kxvef}), 
we have the pull-back 
\begin{align*}
g^* \col (H(X/S,{\cal E}),F)\lo Rg_*((H(X'/S',{\cal E'}),F)
\end{align*} 
fitting into the following commutative diagram 
\begin{equation*}
\begin{CD}
({\cal E}'\otimes_{{\cal O}_{X'}}\Om^{\bul}_{X'/S'},F)@<{\sim}<<
({\cal E}'\otimes_{{\cal O}_{X'}}\Om^{\bul}_{X'/\os{\circ}{S}{}'}[U_{S'}],F) 
@>{\sim}>>(H(X'/S'/\os{\circ}{S}{}',{\cal E}'),F)\\
@V{g^*}VV @V{g^*}VV @V{g^*}VV\\
Rg_*(({\cal E}\otimes_{{\cal O}_X}\Om^{\bul}_{X/S},F))@<{\sim}<<
Rg_*(({\cal E}\otimes_{{\cal O}_X}\Om^{\bul}_{X/\os{\circ}{S}}[U_S],F)) 
@>{\sim}>>Rg_*((H(X/S,{\cal E}),F)). 
\end{CD}
\tag{4.11.2}\label{ali:oobas}
\end{equation*} 
The transitive relation is now obvious. 
\end{proof} 

\begin{coro}
This isomorphism {\rm (\ref{ali:kxef})} 
is contravariantly functorial with respect to $g$ in {\rm (\ref{theo:mra})} 
and the commutative diagrams {\rm (\ref{ali:ooxas})} and {\rm (\ref{cd:kgsx})}.
This functoriality satisfies the transitive relation.  
\end{coro}

\end{document}